\documentclass[11pt]{amsart}
\usepackage{amssymb}

\usepackage{xypic}
\xyoption{all}

\usepackage[linktocpage=true,hypertex]{hyperref}
\newcommand{\myendbibitem}{\relax}

\numberwithin{equation}{section}

\swapnumbers
\newtheorem{thm}[equation]{Theorem}
\newtheorem{prop}[equation]{Proposition}
\newtheorem{lem}[equation]{Lemma}
\newtheorem{cor}[equation]{Corollary}

\theoremstyle{definition}
\newtheorem{defn}[equation]{Definition}
\newtheorem{remark}[equation]{Remark}
\newtheorem{example}[equation]{Example}

\newcommand{\Gr}{\operatorname{Gr}}
\newcommand{\Span}{\operatorname{Span}}
\newcommand{\Spec}{\operatorname{Spec}}
\newcommand{\sdp}{\rtimes}

\newcommand{\bbA}{{\mathbb A}}

\newcommand{\bbZ}{{\mathbb Z}}

\newcommand{\bbQ}{{\mathbb Q}}

\newcommand{\BS}{\operatorname{BS}}

\newcommand{\Gmn}{G_{m,n}}

\newcommand{\inv}{^{-1}}
\newcommand{\Mnm}{(\Mn)^m}
\newcommand{\N}{\operatorname{N}}
\newcommand{\Q}{\operatorname{Q}}
\newcommand{\Z}{\operatorname{Z}}

\newcommand{\GL}{{\operatorname{GL}}}
\newcommand{\SL}{{\operatorname{SL}}}

\newcommand{\PGL}{{\operatorname{PGL}}}
\newcommand{\PGLn}{{\operatorname{PGL}_n}}

\newcommand{\Stab}{\operatorname{Stab}}

\newcommand{\RMaps}{\operatorname{RMaps}}
\newcommand{\lra}{\longrightarrow}
\newcommand{\Mat}{{\operatorname{M}}}
\newcommand{\M}{\operatorname{M}}
\newcommand{\Mn}{\Mat_n}

\newcommand{\trdeg}{\operatorname{trdeg}}
\newcommand{\Sym}{{\operatorname{S}}}
\newcommand{\tr}{\operatorname{tr}}

\newcommand{\GLm}{{\GL_m}}
\newcommand{\numberedpar}{
        \par\refstepcounter{equation}\medskip\noindent{\bf\theequation. }}

\newcommand{\UD}{\operatorname{{\it UD}}}
\newcommand{\UDmn}{\UD(m,n)}

\newcommand{\Zp}{\bbZ/p\bbZ}

\begin{document}

\hbox{}\vspace{-2cm}

\title[Group actions]%
   {Group actions on central simple algebras:\\ a geometric approach }

\date{September 16, 2005} 

\author{Z. Reichstein}
\address{Department of Mathematics, University of British Columbia,
  Vancouver, BC V6T 1Z2, Canada}
\email{reichste@math.ubc.ca}
\urladdr{www.math.ubc.ca/~reichst}
\thanks{Z. Reichstein was supported in part by an NSERC research
  grant.}

\author{N. Vonessen} \address{Department of Mathematical Sciences,
  University of Montana, Missoula, MT 59812-0864, USA}
\email{Nikolaus.Vonessen@umontana.edu}
\urladdr{www.math.umt.edu/vonessen}
\thanks{N.\ Vonessen gratefully acknowledges the support of the
  University of Montana and the hospitality of the University of
  British Columbia during his sabbatical in 2002/2003,
  when part of this research was done.}




\begin{abstract} 
  We study actions of linear algebraic groups on central simple
  algebras using algebro-geometric techniques.  Suppose an algebraic
  group $G$ acts on a central simple algebra $A$ of degree $n$. We are
  interested in questions of the following type: (a) Do the $G$-fixed
  elements form a central simple subalgebra of $A$ of degree $n$?  (b)
  Does $A$ have a $G$-invariant maximal subfield?  (c) Does $A$ have a
  splitting field with a $G$-action, extending the $G$-action on the
  center of $A$?
  
  Somewhat surprisingly, we find that under mild assumptions on $A$
  and the actions, one can answer these questions by using techniques
  from birational invariant theory (i.e., the study of group actions
  on algebraic varieties, up to equivariant birational isomorphisms).
  In fact, group actions on central simple algebras turn out to be
  related to some of the central problems in birational invariant
  theory, such as the existence of sections, stabilizers in general
  position, affine models, etc. In this paper we explain these
  connections and explore them to give partial answers to questions
  (a)---(c).
\end{abstract}

\maketitle

\vspace{-.3in}

\tableofcontents


\section{Introduction}

In this paper we study actions of linear algebraic groups $G$ on
central simple algebras $A$ in characteristic zero. As usual, 
we will denote the center of $A$ by $\Z(A)$ and the subalgebra of
$G$-fixed elements of $A$ by 
\[ A^G = \{ a \in A \, | \, g(a) = a \; \; \forall
g \in G \} \, . \]
We will be interested in questions such as the following:
\begin{quote}%
\begin{itemize} \refstepcounter{equation}\label{q1}
  \item[\llap{(\theequation)\hspace{.25in}}(a)]
   Is $A^G$ a central simple algebra of the same degree as $A$?
  \item[(b)]Does $A$ have a $G$-invariant maximal subfield?
  \item[(c)] Can the $G$-action on $\Z(A)$ be extended to a splitting
    field $L$, and if so, what is the minimal possible value of
    $\trdeg_{\Z(A)} \, L$?
\end{itemize}
\end{quote}
Actions of finite groups on central simple algebras have been
extensively studied in the 1970s and 80s in the context of group
actions on noncommutative rings; for an overview see~\cite{montgomery}.
More recently, torus actions were considered in~\cite{ta} and \cite{RCSA},
and actions of solvable groups in~\cite{vonessen:solv},
all by purely algebraic methods (cf.\ also
\cite{vonessen:memoir,vonessen:trans}).  Inner actions of compact
groups were studied in \cite{sage}.
The purpose of this paper is to introduce a geometric approach
to the subject by relating it to ``birational invariant theory",
i.e., to the study of group actions on algebraic varieties,
up to birational isomorphism.  In particular,
we will see that the questions posed in~\eqref{q1}
are related to some of the central
problems in birational invariant theory, such as
existence of affine models, quotients, stabilizers
in general positions, sections, etc.
(For an overview of birational invariant theory,
see~\cite[Chapters 1, 2, 7]{pv} and \cite[Part 1]{popov2}.)
To make the algebro-geometric techniques applicable, we
always assume that the centers of our simple algebras are
finitely generated field extensions of a fixed algebraically
closed base field $k$ of characteristic zero. All 
algebraic groups are assumed to be linear and defined 
over $k$.

Let $G$ be an algebraic group and $A$ be a finite-dimensional
central simple algebra.  Of course, we are primarily interested 
in studying
$G$-actions on $A$ which respect the structure of $G$ as an algebraic 
(and not just an abstract) group. The following definition is natural 
in the geometric context. 

It is well known that a finitely generated field
extension of $k$ can be interpreted as the field of rational
functions~$k(X)$ on some irreducible variety~$X$, where $X$ is unique
up to birational isomorphism.  Similarly, a central simple algebra~$A$
of degree~$n$ is isomorphic (as a $k$-algebra) to the algebra $k_n(X)$
of $\PGLn$-equivariant rational functions $X \dasharrow\Mn(k)$,
where $X$ is an irreducible variety with a generically free
$\PGLn$-action.  Here $X$ is unique up to birational isomorphism of
$\PGLn$-varieties.  For details, see \cite[Theorem~7.8 and Section~8]{rv4}.  

We will say that a $G$-action 
 on a central simple algebra $A=k_n(X)$ is {\em geometric},
if it is induced by a regular $G$-action on $X$, via
\begin{equation} \label{e1}
(gf)(x)=f(g\inv x)
\end{equation}
for $x \in X$ in general position.
One can check that all rational functions $gf\colon X\dasharrow
\Mn(k)$ lie in $k_n(X)$ (i.e., are $\PGLn$-equivariant) if and only if
the actions of $G$ and $\PGLn$ on $X$ commute.  So a regular
$G$-action on $X$ induces a $G$-action on $A=k_n(X)$ precisely if $X$
is a $G\times\PGLn$-variety.  To sum up:

\begin{defn} \label{def.geometric}
  An action of an algebraic group $G$ on a central simple
  algebra $A$ of degree~$n$ is said to be {\em geometric} if there is
  an irreducible $G \times \PGL_n$-variety $X$ such that $A$ is
  $G$-equivariantly isomorphic to $k_n(X)$.  We will call $X$ {\em the
  associated variety} for this action.
\end{defn}

The second part of the definition makes sense since the associated
variety $X$ is unique up to birational isomorphism (as a $G \times
\PGLn$-variety); see Corollary~\ref{cor:unique}.  Note that 
the $\PGL_n$-action on $X$ is necessarily generically free, 
since $A\simeq k_n(X)$ is a central simple algebra of degree~$n$; see
Lemma~\ref{lem4.1-NEW}. Conversely, any $G \times \PGLn$-variety $X$,
which is $\PGLn$-generically free, is the associated variety for
the geometric action of $G$ on the central simple algebra 
$A = k_n(X)$ given by~\eqref{e1}.

   From an algebraic point of view it is natural to consider 
another class of actions, introduced in~\cite[\S2]{vonessen:solv} 
(and in the special case of torus actions in~\cite[\S5]{RCSA}). 
We shall call such actions {\em algebraic}; for a precise 
definition, see Section~\ref{sect.algebraic}. 
The relationship between algebraic and geometric actions is discussed in
Sections~\ref{sect.algebraic} and \ref{sect.geom=>alg}. In particular,
every algebraic action is geometric; see Theorem~\ref{thm1}.

We are now ready to address the questions posed in~\eqref{q1}, 
in the context of geometric actions.

\begin{thm} \label{thm2}\mbox{}
  Consider a geometric action of an algebraic group $G$ on a
  central simple algebra $A$ of degree~$n$, with associated $G \times
  \PGLn$-variety~$X$.
  \begin{itemize}
  \item[\textup{(a)}]
   The fixed algebra $A^G$ is a central simple
    algebra of degree $n$ if and only if for $x\in X$ in general position,
    \[ \Stab_{G \times \PGLn}(x) \subseteq G \times \{ 1 \}\,.\]

  \item[\textup{(b)}]
   The fixed algebra $A^G$ contains an element with $n$ 
   distinct eigenvalues if and only if for every $x \in X$ 
   in general position there exists a torus $T_x$ of $\PGLn$ such that
    \[ \Stab_{G \times \PGLn}(x) \subseteq G \times T_x\,.\]
\end{itemize}
\end{thm}

We now turn to question~(b) in~\eqref{q1}.
 
\begin{thm} \label{thm4} \mbox{}
  Consider a geometric action of an algebraic group $G$ on a
  central simple algebra $A$ of degree~$n$, with associated $G \times
  \PGLn$-variety~$X$.
  \begin{itemize}
  \item[\textup{(a)}] $A$ has a $G$-invariant maximal \'etale
    subalgebra if and only if there exists a $G \times \PGL_n$-equivariant
    rational map $X \dasharrow \PGL_n/N$, where $N$ is the normalizer
    of a maximal torus in $\PGLn$ and $G$ acts trivially
    on the homogeneous space $\PGLn/N$.

\smallskip
  \item[\textup{(b)}] If $A$ has a $G$-invariant maximal \'etale
    subalgebra, then for every $x \in X$ in general position 
   there exists a maximal torus $T_x$ of $\PGLn$ such that
    \[\Stab_{G \times \PGLn}(x) \subseteq G \times \N(T_x)\,.\]
    Here $\N(T_x)$ denotes the normalizer of $T_x$ in $\PGLn$.

\smallskip
  \item[\textup{(c)}]If the orbit $Gx$ has codimension $< n^2 - n$
    in $X$ for $x \in X$ in general position, then $A$ has no
    $G$-invariant maximal \'etale subalgebras.
\end{itemize}
\end{thm}

Here by an {\em \'etale subalgebra} of $A$ we mean a subalgebra of $A$
which is an \'etale algebra over $\Z(A)$; cf.~\ref{etale:subalgs}.
If~$A$ is a division algebra, the maximal \'etale 
subalgebras are just the maximal subfields.

The converse to Theorem~\ref{thm4}(b) is false in general; 
see Proposition~\ref{prop.example1}.  Note that the points 
of the homogeneous space $\PGLn/N$ parameterize the maximal
tori in $\PGLn$ (see the beginning of~\S\ref{sect.thm4}).  The
converse to part (b) is thus true if and only if the tori $T_x$ can be
chosen so that $x \dasharrow T_x$ is a rational map. We also remark that
Theorem~\ref{thm2}(b) gives a necessary and sufficient condition
for $A$ to have a $G$-invariant maximal \'etale algebra of the form
$\Z(A)[a]$, where $a \in A^G$; see Corollary~\ref{cor.of.thm2}.

Our final result addresses question~(c) in~\eqref{q1}. We begin with 
the following definition.  

\begin{defn} \label{def.G-split}
  Suppose a group $G$ acts on a central simple algebra $A$ of degree
  $n$.  We will say that $A$ is $G$-{\em split}, if $A$ is
  $G$-equivariantly isomorphic to $\Mn(\Z(A)) = \Mn(k) \otimes_k \Z(A)$,
  where $G$ acts via the second factor.  We will say that a
  $G$-equivariant field extension $L/\Z(A)$ is a {\it $G$-splitting
    field} for $A$ if $A \otimes_{\Z(A)} L$ is $G$-split.
\end{defn}

Note that if $G$ acts trivially on $A$, then a $G$-splitting field is
just a splitting field for $A$ in the usual sense.  Note also that
a $G$-action on a split central simple algebra (i.e., a matrix algebra
over a field) need not be $G$-split (cf.\ Example~\ref{ex:m2}).

\begin{thm} \label{thm3}
  Every geometric action of an algebraic group $G$ on a central
  simple algebra $A$ of degree $n$ has a $G$-splitting field 
  of the form $L = k(X_0)$, where $X_0$ is a $G$-variety 
  and $\trdeg_{\Z(A)}(L) = n^2 -1$.  Moreover, if $G$ acts 
  algebraically on $A$, then $X_0$ can, in addition, be chosen 
  to be affine.
\end{thm}

In general, the value of $\trdeg_{\Z(A)} L$ given in
Theorem~\ref{thm3} is the smallest possible; see
Proposition~\ref{prop.split}(b). If $G$ is connected, we  
give a different construction of $G$-splitting fields in 
Section~\ref{sect3b}.

At the end of the paper we will present four examples illustrating our
main results, Theorems~\ref{thm2}, \ref{thm4}, and \ref{thm3}, and two
appendices.  Appendix~$A$ deals with inner actions on division
algebras which need not be geometric, while Appendix~B treats regular
actions of algebraic groups (see Definition~\ref{def.regular}) on
prime affine PI-algebras.  Using Theorem~\ref{thm3}, we show that such
actions are ``induced'' by regular actions on commutative domains.
Further results on geometric actions will appear in the
paper~\cite{rv5}.

\section{Preliminaries}
\label{sect.prel}

\numberedpar \emph{Conventions.}  We work over a fixed algebraically
closed base field~$k$ of characteristic zero. All algebras are
$k$-algebras, and division algebras and central simple algebras are
assumed to be finite-dimensional over their centers, which in turn are
assumed to be finitely generated field extensions of~$k$.  All actions
on algebras are by $k$-algebra automorphisms.  Algebraic groups are
always assumed to be linear algebraic groups over~$k$, and $G$ will always
denote an algebraic group. 
Regular actions are meant to be regular
over~$k$; similarly for algebraic actions (see
Definition~\ref{defn:algebraic-action}).
If $K$ is a field, we shall denote
the algebra of $n \times n$ matrices over $K$ by $\Mn(K)$.
If $K = k$, we will write $\Mn$ in place of $\Mn(k)$. We will
sometimes view $\Mn$ as a $k$-algebra and sometimes as 
an algebraic variety, isomorphic to the affine space $\bbA^{n^2}$.

\numberedpar \emph{$G$-varieties.}  By a $G$-variety $X$ we mean an
algebraic variety with a regular action of $G$.  By a morphism $X \lra
Y$ of $G$-varieties, we mean a $G$-equivariant morphism. The notions
of isomorphism, rational map, birational isomorphism, etc.\ of
$G$-varieties are defined in a similar manner.  As usual, given a
$G$-action on $X$, we denote the orbit of $x \in X$ by $Gx$ and the
stabilizer subgroup of $x$ by $\Stab_G(x) \subseteq G$. Throughout this paper
we use~\cite{pv} as a reference for standard notions from invariant
theory, such as rational and categorical quotients, stabilizers in
general position, sections, etc.

\begin{defn} \label{def.invariant}
We shall say that a $G$-action on $X$ is 

\smallskip
(a) {\em faithful} if every $1 \ne g \in G$ acts nontrivially on $X$,

\smallskip
(b) {\em generically free} if $\Stab_G(x) = \{ 1 \}$ for $x \in X$
in general position, and

\smallskip
(c) {\em stable} if the orbit $G x$ is closed in $X$ 
for $x \in X$ in general position.
\end{defn}

\begin{lem} \label{lem.diag}
  Suppose the group $G$ is either {\upshape(a)} finite or
  {\upshape(b)} diagonalizable.  Then every faithful irreducible
  $G$-variety $X$ is generically free.
\end{lem}

\begin{proof} (a) Since the $G$-action is faithful,
  $X^g = \{ x \in X \, | \, gx = x \} \neq X$ for every $1\neq g \in G$.
Since each $X^g$ is a closed subvariety of $X$, every point of
the Zariski dense open subset $X - \cup_{1\neq g \in G} \, X^g$
has a trivial stabilizer in $G$.

\smallskip
Part (b) is an immediate corollary of a theorem of
Richardson~\cite[Theorem 9.3.1]{richardson};
see also~\cite[Theorem 7.1]{pv}.
\end{proof}

The following example shows that, contrary to the assertion
in~\cite[Proposition 7.2]{pv}, Lemma~\ref{lem.diag} fails
if we only assume that the connected component of $G$
is a torus. We shall return to this example in~\S\ref{sect.ex2x2}.

\begin{example} \label{ex.O_2}
  Consider the natural linear action of the orthogonal group $G = O_2$
  on $\bbA^2$. This action is faithful but not generically free:
  $\Stab_G(v)$ has order 2 for $v \in k^2$ in general position.
  Indeed, for every non-isotropic vector~$v$ in $k^2$, there is a
  unique non-trivial element of $O_2$, leaving $v$ invariant; this
  element is the orthogonal reflection in $v$.  Note also $O_2$ is a
  semidirect product of a one-dimensional torus with $\bbZ/2 \bbZ$.
\end{example}

\begin{lem}[Popov] \label{lem.popov}
  Let $G$ be a reductive group, $X$ be an affine $G$-variety and $V$
  be a $G$-representation. Suppose the $G$-orbit of $x \in X$ is
  closed in $X$ and $\Stab(x) \subseteq \Stab(v)$ for some $v \in V$.
  Then there exists a $G$-invariant morphism $f \colon X \lra V$ such
  that $f(x) = v$.
\end{lem}

\begin{proof}
  In the case where $\Stab(x) = \{ 1 \}$, this lemma is stated and
  proved in \cite[Theorem 1.7.12]{popov2}.  The same argument goes
  through in our slightly more general setting.
\end{proof}

\numberedpar\label{maps} \emph{Algebras of rational maps.}  If $X$ is
a $\PGLn$-variety, we will denote by $\RMaps_{\PGLn}(X, \Mn)$ the
$k$-algebra of $\PGLn$-equivariant rational maps $f \colon X
\dasharrow \Mn$, with addition and multiplication induced from $\Mn$.

\begin{lem} \label{lem4.1-NEW}
  Let $Y$ be an irreducible $\PGLn$-variety.  Then the following are
  equivalent:
\begin{itemize}
\item[\textup(a)]The $\PGLn$-action on $Y$ is generically free.
\item[\textup(b)]$A=\RMaps_{\PGLn}(Y,\Mn)$ is a central simple algebra
  of degree~$n$.
\end{itemize}
If {\upshape(a)} and {\upshape(b)} hold then the center of $A$ is $\RMaps_{\PGLn}(Y, k) =
k(Y)^{\PGLn}$. Here elements of $k$ are identified with scalar matrices
in $\Mn$.
\end{lem}

\begin{proof}
  (b) $\Rightarrow$ (a):
  Note that the center of $A$ contains $k(Y)^{\PGLn}$.
  Choose $f_1, \ldots, f_{n^2}\in A$ which are
  linearly independent over $k(Y)^\PGLn$.  By \cite[Lemma
  7.4]{reichstein}, $f_1(y),\ldots,f_{n^2}(y)$ are $k$-linearly
  independent in $\Mn$ for $y\in Y$ in general position.  Now consider
  the $\PGLn$-equivariant rational map
  \[ f=(f_1,\ldots,f_{n^2}) \colon Y \dasharrow (\Mn)^{n^2} \, .\]
  For $y \in Y$ in general position, $\Stab(f(y)) = \{ 1 \}$, so
  that also $\Stab(y) = \{ 1 \}$.  Hence $Y$ is $\PGLn$-generically
  free.

  The implication (a) $\Rightarrow$ (b) and the last assertion
  of the lemma are proved in~\cite[Lemma 8.5]{reichstein} (see also
  \cite[Definition 7.3 and Lemma 9.1]{reichstein}).
\end{proof}

\smallskip
If the $\PGLn$-action on $X$ is generically free, we will denote 
the central simple algebra $\RMaps_{\PGLn}(X, \Mn)$ by $k_n(X)$.

\numberedpar{\em Maximal \'etale subalgebras.}
\label{etale:subalgs}
Let $A$ be a central simple algebra of degree~$n$.  By an \'etale
subalgebra of $A$ we mean a subalgebra of $A$ which is an \'etale
algebra over $\Z(A)$, i.e., a finite direct sum of (separable) field
extensions of $\Z(A)$. Note that since we are working in
characteristic zero, the term ``\'etale" could be replaced by
``commutative semisimple''. We are interested in maximal \'etale
subalgebras, i.e., \'etale subalgebras $E$ of $A$ satisfying the
following equivalent conditions:
\begin{itemize}
\item[(a)]$\dim_{\Z(A)} \, E = \deg(A)$,
\item[(b)]$E$ is maximal among commutative subalgebras of $A$;
\end{itemize}
cf.~\cite[Exercise 7.1.1]{rowen:ringII}.  Using the double centralizer
theorem, one easily verifies that every \'etale subalgebra of $A$ is
contained in a maximal \'etale subalgebra, see, e.g.,
\cite[Theorem~4.10 and Exercise~4.6.12]{jacobson:BAII} and
\cite[Exercise 7.1.2]{rowen:ringII}.  Of course, if $A$ is a division
algebra, then maximal \'etale subalgebras are just maximal subfields.

We will repeatedly use the following characterization of maximal
\'etale subalgebras, which follows easily from
\cite[\S V.7.2, Proposition 3]{bourbaki:algebra}.

\begin{lem}\label{lem:etale}
  Let $A$ be a central simple algebra of degree~$n$ with center $K$.
  Let $a\in A$. Then $K[a]$ is a maximal \'etale subalgebra of $A$ if
  and only if the eigenvalues of $a$ are distinct.
  \qed
\end{lem}

\section{The uniqueness of the associated variety}
\label{sect.uniqueness}

Recall that given a generically free $\PGLn$-variety $X$, we write $A
= k_n(X)$ for the algebra of $\PGLn$-equivariant functions $a \colon X
\dasharrow \Mn$. A $\PGLn$-equivariant dominant rational map $f \colon
X' \dasharrow X$ induces an embedding $f^* \colon A \hookrightarrow
A'$ of central simple algebras, where $A' = k_n(X')$ and $f(a) = a
\circ f \colon X' \dasharrow \Mn$.

We now deduce a simple consequence of the functoriality of the maps
$X \mapsto k_n(X)$ and $f \mapsto f^*$; see~\cite[Theorem 1.2]{rv4}. 
Recall that if $X$ has a $G$-action, which commutes with 
the $\PGLn$-action, then~\eqref{e1} defines a $G$-action
on $A = k_n(X)$, which we call geometric.

\begin{lem} \label{lem:unique}
  Let $X$ and $X'$ be $G \times \PGLn$-varieties, which are
  $\PGLn$-generically free.
\begin{itemize}
\item[\upshape(a)]If $f \colon X' \dasharrow X$ is a dominant rational
  map of $G \times \PGLn$-varieties then the induced embedding $f^*
  \colon k_n(X) \hookrightarrow k_n(X')$ of central simple algebras is
  $G$-equivariant.
\item[\upshape(b)]Every $G$-equivariant embedding $j \colon k_n(X)
  \hookrightarrow k_n(X')$ induces a dominant rational $G \times
  \PGLn$-equivariant map $j_* \colon X' \dasharrow X$.
\end{itemize}
\end{lem}

\begin{proof} (a) By~\cite[Theorem 1.2]{rv4}, since the diagram
  \[ \xymatrix{
    X' \ar@{-->}[r]^{f} \ar@{->}[d]^{g} & X \ar@{->}[d]^{g} \cr
    X' \ar@{-->}[r]^{f}                 & X } \]
commutes for every $g \in G$, so does the induced diagram
  \[ \xymatrix{
    k_n(X) \ar@{->}[r]^{f^*} \ar@{->}[d]^{g\inv} 
                           & k_n(X') \ar@{->}[d]^{g\inv} \cr 
    k_n(X) \ar@{->}[r]^{f^*} & k_n(X'). } \]

(b) Conversely, since the diagram
  \[ \xymatrix{
    k_n(X) \ar@{->}[r]^{j} \ar@{->}[d]^{g\inv} 
                           & k_n(X') \ar@{->}[d]^{g\inv} \cr 
    k_n(X) \ar@{->}[r]^{j} & k_n(X') } \]
commutes, so does the induced diagram
\[    \xymatrix{
    X' \ar@{-->}[r]^{j_*} \ar@{->}[d]^{g} & X \ar@{->}[d]^{g} \cr
    X' \ar@{-->}[r]^{j_*}                 & X } \]
\end{proof}

\begin{cor} \label{cor:unique}
Given a geometric action of an algebraic group $G$ on a
  central simple algebra $A$, the $G \times \PGLn$-variety associated
  to this action is unique up to birational isomorphism.
\end{cor}

\begin{proof}
  Suppose two $G \times \PGLn$-varieties $X$ and $X'$ are both
  associated varieties for this action, i.e., $k_n(X)$ and $k_n(X')$
  are both $G$-equi\-var\-i\-antly isomorphic to $A$. In other words,
  there are mutually inverse $G$-equivariant algebra isomorphism
  $i \colon k_n(X) \stackrel{\simeq}{\lra} k_n(X')$ and
  $j \colon k_n(X') \stackrel{\simeq}{\lra} k_n(X)$.
  Applying Lemma~\ref{lem:unique}, $i$ and $j$ induce mutually inverse
  dominant $G \times \PGLn$-equivariant rational map
  $i_* \colon X' \stackrel{\simeq}{\dasharrow} X$ and
  $j_* \colon X \stackrel{\simeq}{\dasharrow} X'$.
  We conclude that $X$ and $X'$ are birationally isomorphic
  $G \times \PGLn$-varieties.
\end{proof}

\begin{example} \label{ex.conjugation}
  Let $G$ be a subgroup of $\PGL_n$, and consider the conjugation
  action of $G$ on $A = \Mn(k)$. We claim that the associated $G
  \times \PGLn$-variety for this action is $X = \PGLn$, with $G$
  acting by translations on the right and $\PGLn$ acting by
  translations on the left.  More precisely, for $(g,h)\in
  G\times\PGLn$ and $x\in X$, $(g,h)\cdot x=hxg\inv$.  Consequently
  for $f\in k_n(X)$,
  \[(g\cdot f)(x)=f\bigl((g,1)\inv\cdot x\bigr)=f(xg)\,,\]
  see \eqref{e1}.  Note that since $X$ is a single $\PGLn$-orbit,
  every $\PGLn$-equivariant rational map $f \colon \PGLn \dasharrow
  \Mn$ is necessarily regular.  It is
  now easy to check that the $k$-algebra isomorphism
  \[ \phi \colon k_n(X) 
       = \RMaps_{\PGLn}(\PGLn, \Mn) \stackrel{\simeq}{\lra} A=\Mn \, \]
  given by $\phi(f) = f(1)$ is $G$-equivariant.
\end{example}

\begin{example} \label{ex.universal}
  Let $m\geq 2$, and consider the $\PGLn$-variety $X = \Mnm$, where
  $\PGLn$ acts by simultaneous conjugation, i.e., via
  \[g \cdot (a_1, \dots, a_m) = (g a_1 g\inv, \dots, g a_m g\inv)\, . \] 
  Since $m \ge 2$, this action is generically free. The associated
  division algebra $k_n(X)$ is called the {\em universal division
    algebra} of $m$ generic $n \times n$-matrices and is denoted by
  $\UDmn$.  Identify the function field of $X$ with $k(x_{ij}^{(h)})$,
  where for each $h = 1, \dots, m$, $x_{ij}^{(h)}$ are the $n^2$
  coordinate functions on copy number~$h$ of $\Mn$, and identify the
  algebra of all rational maps $X \dasharrow \Mn$ with
  $\Mn\bigl(k(x_{ij}^{(h)})\bigr)$. Now we can think of $\UDmn$ as the
  division subalgebra of $\Mn\bigl(k(x_{ij}^{(h)})\bigr)$ generated by
  the $m$ {\em generic $n \times n$ matrices} $X^{(h)} =
  (x_{ij}^{(h)})$, $h = 1, \dots, m$.  Here $X^{(h)}$ corresponds to
  projection $\Mnm \lra \Mn$ given by $(a_1, \dots, a_m) \mapsto a_h$.
  For details of this construction, see~\cite[Section 2]{procesi2}
  or~\cite[Theorem 5]{lb}.

  Now observe that the $\GL_m$-action on $X = \Mnm$ given by
  \begin{equation}\label{eq.action.on.Mnm}
  g \cdot (a_1, \dots, a_m) = (\sum_{j = 1}^m g_{1j} a_j, \dots,
  \sum_{j = 1}^m g_{mj} a_j) 
  \end{equation}
  commutes with the above
  $\PGLn$-action. Here $g = (g_{ij}) \in \GL_m$, with $g_{ij} \in k$.
  Using formula~\eqref{e1}, we see that this gives rise to a
  $\GL_m$-action on $\UDmn$ such that for $g\in\GL_m$,
  \begin{equation}\label{eq.action.on.UDmn}
    g \cdot X^{(h)}=\sum_{j = 1}^m g'_{hj} X^{(j)}\,,
  \end{equation}
  where $g^{-1} = (g'_{ij})$.  In other words, this $\GL_m$-action on
  $\UDmn$ is geometric, with associated $G \times \PGLn$-variety $X =
  \Mnm$.  We will return to this important example later in this paper
  (in Example~\ref{ex.universal.alg} and Sections~\ref{sect.ex.UDmn}
  and~\ref{sect.ex2x2}), as well as in \cite{rv5}.
\end{example}

\begin{remark}\label{rem.Gmn}
  The $k$-subalgebra of $\UDmn$ generated by $X^{(1)},\ldots,X^{(m)}$
  is called the {\em generic matrix ring} generated by $m$ generic
  $n\times n$ matrices; we denote it by $\Gmn$.  Note that the
  action~\eqref{eq.action.on.UDmn} of $\GLm$ on $\UDmn$ restricts to
  an action on $\Gmn$. Consequently, the $\GLm$-action on $\Gmn$
  is induced by the $\GLm$-action on $\Mnm$ in the sense of
  formula~\eqref{e1}.
\end{remark}

\section{Brauer-Severi varieties}

Let $A/K$ be a central simple algebra of degree~$n$. Throughout much
of this paper, 
we associate to $A$ a $\PGLn$-variety $X/k$ such that $A$ is the algebra
of $\PGLn$-equivariant rational maps $X \dasharrow \Mn(k)$. Another 
variety that can be naturally associated to $A$ is the Brauer-Severi 
variety $\BS(A)$, defined over $K$. Any algebra automorphism 
$g \colon A \lra A$, defined over the base field $k$, 
induces $k$-automorphisms of $K$ and $\BS(A)$ such that the diagram
  \[ \xymatrix{
    \BS(A) \ar@{->}[r]^{g_*} \ar@{->}[d] &  \BS(A) \ar@{->}[d] \cr
    \Spec(K) \ar@{->}[r]^{(g_{| \, K})_*} & \Spec(K) \, ,  } \]
commutes; conversely, $g$ can be uniquely recovered from 
this diagram.  If a group $G$ acts on $A$, it is natural to ask if
$\BS(A)$ can be $G$-equivariantly represented by an algebraic variety
over $k$. In this short section we will address this question, 
following a suggestion of the referee. Our main result,
Proposition~\ref{prop.bs} below, will not be used in the sequel.  

\begin{prop} \label{prop.bs} Consider a geometric action $\phi$ of 
an algebraic group $G$ on a central simple algebra $A/K$ of degree~$n$.
Then there exists a morphism $\sigma \colon S \lra Y$ 
of irreducible $G$-varieties \textup{(}of finite type over
$k$\textup{)} such that 
\begin{itemize}
\item[(a)] $S$ is a Brauer-Severi variety over $Y$;
\item[(b)] $k(Y) = K$ and $\sigma^{-1}(\eta)$ is the Brauer-Severi
  variety of $A$, where $\eta$ is the generic point of $Y$;
\item[(c)] the $G$-actions on $S$ and $Y$ induce the action $\phi$ on $A$.
\end{itemize}
\end{prop}

\begin{proof} Let $X$ be the $G \times \PGLn$-variety associated to
$\phi$ and $H$ be the maximal parabolic
subgroup of $\PGLn$ consisting of matrices of the form 
\[ \begin{pmatrix} * & 0 & \dots & 0 \\
* & * & \dots & * \\
\vdots & \vdots & & \vdots \\
* & * & \dots & * \end{pmatrix} \, . \]
Consider the natural dominant
rational map $\sigma \colon X/H \dasharrow X/\PGLn$ given by
the inclusion $k(X)^{\PGLn} \hookrightarrow k(X)^H$.
Recall that the rational quotient varieties $X/H$ and $X/\PGL$
are a priori only defined up to birational isomorphism. However, 
we can choose models for these varieties such that
the induced $G$-actions are regular; 
cf.~\cite[Proposition 2.6 and Corollary 1.1]{pv}. 
For notational convenience, we will continue to denote these 
$G$-varieties by $X/H$ and $X/\PGLn$.
Note also that since the actions of $G$ and $\PGLn$ 
on $X$ commute, the resulting map 
$\sigma \colon X/H \dasharrow X/\PGLn$ is $G$-equivariant. 

By~\cite[Section 9]{rv4}, $X/H$ is a Brauer-Severi 
variety over a dense open subset $U$ of $X/\PGLn$, and is isomorphic
to $\BS(A)$ over the generic point of $X/\PGLn$. Since $\sigma$
is $G$-equivariant, $X/H$ is a Brauer-Severi
variety over $g(U)$, for every $g \in G$.
Setting $Y$ to be the union of the $g(U)$ inside $X/\PGLn$, 
as $g$ ranges over $G$, and setting $S$ to be the preimage 
of this set in $X/H$, we obtain a $G$-equivariant 
morphism $\sigma \colon S \lra Y$ with desired properties.
\end{proof}

\section{Algebraic actions}
\label{sect.algebraic}

\begin{defn} \label{def.regular}
We shall say that the action of an algebraic group $G$ on a (not
necessarily commutative) $k$-algebra $R$ is {\em regular}
\footnote{Such actions are usually called rational; we prefer the term
  regular, since the term ``rational action" has a different meaning
  in the context of birational invariant theory.}, if every
finite-dimensional $k$-subspace of $R$ is contained in a $G$-invariant
finite-dimensional $k$-subspace $V$, such that the $G$-action on $V$
induces a homomorphism $G \lra \GL(V)$ of algebraic groups.  
\end{defn}

Every regular action of a connected algebraic group on a division
algebra (or even a field) must be trivial (see, e.g.,
\cite[A.1]{vonessen:trans}), so this notion is too restrictive for our
purposes. However, it naturally leads to the following definition,
made in~\cite[\S2]{vonessen:solv}.  (The special case where $G$ is a
torus had been considered earlier in~\cite[\S5]{RCSA}.)
 
\begin{defn}\label{defn:algebraic-action}
  Let $G$ be an algebraic group acting on a $k$-algebra $A$ by
  $k$-algebra automorphisms.  We call the action {\em algebraic}
  \footnote{In \cite{vonessen:solv}, $S$ is not required to be
    central; it is, however, proved there that $S$ can always be
    chosen to be central if $A$ is a central simple algebra.}  (over
  $k$) if there is a $G$-invariant subalgebra $R$ of $A$ and a
  $G$-invariant multiplicatively closed subset $S$ of $R$ consisting
  of central nonzerodivisors of $R$ such that (1)~$G$ acts regularly
  on~$R$, and (2)~$A=RS\inv$.
\end{defn}

Note that a regular action on $A$ is algebraic (use $S=\{1\}$). 
We shall be primarily interested in the case where $A$ is a central 
simple algebra; in this case $R$ is an order in $A$ (and in particular,
$R$ is prime). For basic properties of algebraic actions, 
see~\cite[\S2]{vonessen:solv}.

The purpose of this section is to investigate the relationship
between algebraic and geometric actions 
(cf.\ Definition~\ref{def.geometric}).

\begin{thm}\label{thm1} 
{\upshape (a)} Algebraic actions are geometric.

\smallskip {\upshape(b)} Let $G$ be an algebraic group acting
geometrically on a central simple algebra~$A$ of degree~$n$.  Then the
action of $G$ on $A$ is algebraic if and only if there is an
associated $G\times\PGLn$-variety $X$ with the following two
properties: $X$ is affine, and the $\PGLn$-action on $X$ is stable
{\upshape(}cf.\ Definition~{\upshape\ref{def.invariant}(c))}.
\end{thm}

We begin with a result which is a $G$-equivariant version
of~\cite[Theorem~6.4]{rv4}. 

\begin{prop}\label{thm1-new-prop}
  Let $G$ be an algebraic group acting regularly on a finite\-ly
  generated prime $k$-algebra~$R$ of PI-degree~$n$.  Then there is an
  $n$-variety $Y$ with a regular $G$-action such that $R$ is
  $G$-equivariantly isomorphic to $k_n[Y]$.
\end{prop}

See~\cite[3.1]{rv4} for the definition of $k_n[Y]$, the PI-coordinate
ring of $Y$.  The action of $G$ on $k_n[Y]$ is induced from the action
of $G$ on $Y$ as in formula~\eqref{e1}.

\begin{proof}
  We may assume that $G$ acts faithfully on $R$.  There is a
  finite-dimensional $G$-stable $k$-subspace $W$ of $R$ which
  generates $R$ as a $k$-algebra.  Set $m=\dim_k(W)$, and consider the
  generic matrix ring $\Gmn$ with its $\GLm$-action as in
  Remark~\ref{rem.Gmn}.  Denote by $V$ the $k$-subspace of $\Gmn$
  generated by the $m$ generic $n\times n$ matrices.  Let
  $\psi_0\colon V\to W$ be a $k$-vector space isomorphism.  Define a
  regular action of $G$ on $V$ by making $\psi_0$ $G$-equivariant.
  The action of $G$ on $V$ extends to a regular action on $\Gmn$.  By
  the universal mapping property of $\Gmn$, $\psi_0$ extends to a
  $G$-equivariant surjective $k$-algebra homomorphism $\psi\colon\Gmn
  \to R$.  Replacing $G$ by an isomorphic subgroup of $\GLm,$ we may
  assume that $G$ acts on $V$ as in \eqref{eq.action.on.UDmn}.  Then
  the action of $G$ on $\Gmn$ is induced (as in~\eqref{e1}) from the
  action of $G$ on $\Mnm$ given by \eqref{eq.action.on.Mnm}.  Note
  that the actions of $G$ and $\PGLn$ on $\Mnm$ commute.
  
  Let $I$ be the kernel of $\psi$, and let $Y={\mathcal
    Z}(I)\subset\Mnm$ be the irreducible $n$-variety associated to
  $I$, see \cite[Corollary 4.3]{rv4}.  Note that $Y$ is $G$-stable for
  the action of $G$ on $\Mnm$.  By~\cite[Proposition~5.3]{rv4},
  ${\mathcal I}(Y)=I$, so that $R$ is $G$-equivariantly isomorphic to
  $k_n[Y]=\Gmn/{\mathcal I}(Y)=\Gmn/I$.
\end{proof}

\begin{proof}[Proof of Theorem~\ref{thm1}]
  (a) Let $G$ be an algebraic group acting algebraically on a central
  simple algebra~$A$ of degree~$n$.  Let $R$ be a $G$-stable finitely
  generated prime PI-algebra contained in $A$ such that $A$ is the
  total ring of fractions of $R$.  By Proposition~\ref{thm1-new-prop},
  there is an $n$-variety $Y$ with a regular action of $G$ such that
  $R$ is $G$-equivariantly isomorphic to $k_n[Y]$.  Then $A$ is
  $G$-equivariantly isomorphic to the total ring of fractions of
  $k_n[Y]$, i.e., to $k_n(Y)$, see~\cite[Proposition~7.3]{rv4}.  As
  the proof of Proposition~\ref{thm1-new-prop} showed, $\Mnm$ is a
  $G\times\PGLn$-variety (where $G$ acts via some subgroup of $\GLm$
  as in \eqref{eq.action.on.Mnm}), and $Y$ is a $G$-stable subset of
  $\Mnm$.  Hence the closure $X$ of $Y$ in $\Mnm$ is an affine
  $G\times \PGLn$-variety.  It is clear that $k_n(Y)$ and $k_n(X)$ are
  $G$-equivariantly isomorphic, and that the $\PGLn$-action on $X$ is
  generically free and stable.  So $G$ acts geometrically on $A$, and
  the associated $G\times\PGLn$-variety $X$ has the two additional
  properties from part~(b).

\smallskip

(b) If the action of $G$ on $A$ is algebraic then an 
  associated $G\times\PGLn$-variety $X$ with desired properties
  was constructed in the proof of part (a).
  
  Conversely, assume that there is an associated
  $G\times\PGLn$-variety $X$ which is affine and on which the
  $\PGLn$-action is stable.  We may assume that $A=k_n(X)$.  So $A$ is
  a central simple algebra with center $K=k(X)^{\PGLn}$; cf.\ 
  Lemma~\ref{lem4.1-NEW}.  Since $X$ is affine, and since
  $\PGLn$-orbits in $X$ in general position are closed, $k[X]^{\PGLn}$
  separates $\PGLn$-orbits in general position, so that
  $\Q(k[X]^{\PGLn})=k(X)^{\PGLn}=K$; see~\cite[Lemma~2.1]{pv}. (Here
  $\Q$ stands for the quotient field.) Denote by $R$ the
  subalgebra of $A$ consisting of the regular $\PGLn$-equivariant maps
  $X\to\Mn$.  It is clearly $G$-invariant.  Note that $G$ acts
  regularly on $k[X]$.  Consequently, $G$ acts regularly on
  $\Mn(k[X])$, the set of regular maps $X\to\Mn$.  Hence, $G$ also
  acts regularly on its subalgebra $R$.  It remains to show that $R$
  is a prime subalgebra of $A$, and that its total ring of fractions
  is equal to $A$.
  
  Let $v\in(\Mn)^2$ be a pair of matrices generating $\Mn$ as
  $k$-algebra, and let $x\in X$ be such that its stabilizer in $\PGLn$
  is trivial and such that its $\PGLn$-orbit is closed.  Then by
  Lemma~\ref{lem.popov}, there is a $\PGLn$-equivariant regular map
  $X\to(\Mn)^2$ such that $f(x)=v$.  Write $f=(f_1,f_2)$, where $f_1$
  and $f_2$ are $\PGLn$-equivariant regular maps $X \lra \Mn$, i.e.,
  elements of $R$.  Since $f_1(x)$ and $f_2(x)$ generate $\Mn$, the
  central polynomial~$g_n$ (\cite[p.~26]{rowen:PI}) does not vanish
  on~$R$.  Since $g_n$ is $t^2$-normal, it vanishes on every proper
  $K$-subspace of $A$, see \cite[1.1.35]{rowen:PI}.  Consequently
  $RK=A$, and $R$ is prime and has PI-degree~$n$.  Clearly, $R$
  contains $k[X]^{\PGLn}$.  Since $R\Q(k[X]^{\PGLn})=RK=A$, $\Q(R)=A$.
  Hence, $G$ acts algebraically on~$A$.
\end{proof}

\begin{example}\label{ex.universal.alg}
  It follows easily from Definition~\ref{defn:algebraic-action} that
  the action~\eqref{eq.action.on.UDmn} of $\GL_m$ on $\UDmn$ is
  algebraic.  So by Theorem~\ref{thm1}(b), there is an associated
  $\GL_m\times\PGLn$-variety $X$ with the following two properties:
  $X$ is affine, and the $\PGLn$-action on $X$ is stable.  Indeed, the
  natural associated variety $X=\Mnm$ has these properties.
\end{example}
 
\section{Proof of Theorem~\ref{thm2}}

We begin with the following simple observation:

\begin{remark}\label{rem:first-results}
  Consider a geometric action of an algebraic group $G$ on a central
  simple algebra $A$, with associated $G \times \PGLn$-variety $X$.
  Elements of $A$ are thus $\PGLn$-equivariant rational maps
  $a \colon X \dasharrow \M_n$. Such an element is 
  $G$-fixed if and only if it factors through the rational
  quotient map $X \dasharrow X/G$.  In other words, $A^G$ is
  isomorphic to $\RMaps_{\PGL_n}(X/G, \Mn)$.
\end{remark}

We are now ready to proceed with the proof of Theorem~\ref{thm2}.

\smallskip
(a) We may assume that $A=k_n(X)$. Combining 
  Remark~\ref{rem:first-results} with Lemma~\ref{lem4.1-NEW},
  we see that $A^G$ is a central simple algebra of degree $n$ if
  and only if $Y = X/G$ is a generically free $\PGL_n$-variety. The 
  latter condition is equivalent to 
  $\Stab_{G \times \PGLn}(x)\subseteq G \times \{ 1 \}$
  for $x \in X$ in general position.

\smallskip

(b)  First suppose that there is an $a\in A^G$ with $n$ distinct eigenvalues.
  Adding to $a$ some constant in $k$, we may assume that
  the eigenvalues of $a$ are distinct and nonzero.  Hence for $x\in X$
  in general position, the eigenvalues of $a(x)\in\Mn$ are also
  distinct and nonzero.  The stabilizer of $a(x)$ in $\PGLn$ is thus a
  maximal torus $T_x$ of $\PGLn$.  Let $(g, p) \in \Stab_{G \times
    \PGLn}(x)$.  Then $a(x)=g(a)(x) = a(g^{-1}(x)) = a(p(x)) = p\, a(x)
  p^{-1}$.  Thus $p\in T_x$, so that $\Stab_{G \times \PGLn}(x)
  \subseteq G \times T_x$.

  \smallskip
  
  We will now prove the converse. Assume $\Stab_{G \times \PGLn}(x)$
  is contained in $G \times T_x$ for some torus $T_x$ of $\PGLn$
  (depending on $x$).  Denote by $Y$ the rational quotient
  $\PGLn$-variety $X/G$.  To produce an $a\in A^G$ with distinct
  eigenvalues, it suffices to construct a $\PGLn$-equivariant rational
  map $a \colon Y \dasharrow \Mn$ whose image contains a matrix with
  distinct eigenvalues.  By our assumption, $\Stab_{\PGLn}(y)$ is
  contained in a torus $T_x \subset \PGL_n$ for $y \in Y$ in general
  position. Hence, $\Stab_{\PGLn}(y)$ is diagonalizable (and, in
  particular, reductive).  By \cite[Theorem 1.1]{rv3}, after replacing
  $Y$ by a birationally equivalent $\PGLn$-variety, we may assume that
  $Y$ is affine and the $\PGLn$-action on $Y$ is stable.
  
  We are now ready to construct a map $a \colon Y \dasharrow \Mn$ with
  the desired properties. Let $y \in Y$ be a point whose orbit is
  closed and whose stabilizer $S$ is diagonalizable, and let $v
  \in\Mn$ be a matrix with distinct eigenvalues. Then $\Stab(v)$ is a
  maximal torus in $\PGLn$; after replacing $v$ by a suitable
  conjugate, we may assume $S \subseteq \Stab(v)$. Now
  Lemma~\ref{lem.popov} asserts that there exists a
  $\PGLn$-equivariant morphism $a \colon Y \lra \Mn$ such that $a(y) =
  v$. This completes the proof of Theorem~\ref{thm2}.  \qed

\begin{example}\label{ex:m2}
Let $G$ be a subgroup of $\PGL_n$, acting by conjugation
on $A = \Mn(k)$. The associated
variety for this action is $X = \PGLn$, with $G \times \PGLn$ 
acting on it by $(g, h) \cdot x = h x g^{-1}$; see
Example~\ref{ex.conjugation}.
Since all of $X$ is a single $\PGLn$-orbit, the stabilizer of any $x \in X$
is conjugate to the stabilizer of $1_{\PGLn}$, which is easily seen to be
$\{ (g, g) \, | \, g \in G \}$.
So in this setting, Theorem~\ref{thm2}(b) reduces 
to the following familiar facts:

\begin{itemize}
  \item[(a)]$\Mn(k)^G = \Mn(k)$ if and only if $G = \{ 1 \}$, and
    
  \item[(b)]$\Mn(k)^G$ contains an element with $n$ distinct
    eigenvalues if and only if $G$ centralizes a maximal torus in
    $\GL_n$, i.e., if and only if $G$ is contained in maximal torus of
    $\PGLn$.
\end{itemize}
\end{example}

Using Lemma~\ref{lem:etale}, we can  
rephrase Theorem~\ref{thm2}(b) in a way that makes its 
relationship to Question~\ref{q1}(b) more transparent.

\begin{cor}\label{cor.of.thm2}
  Consider a geometric action of an algebraic group $G$ on a
  central simple algebra $A$ of degree~$n$, with associated $G \times
  \PGLn$-variety~$X$.  The following conditions are equivalent.
  \begin{itemize}
  \item[(a)]$A$ has a maximal \'etale
    subalgebra~$E$ of the form $E=\Z(A)[a]$ for some $a\in A^G$.
  \item[(b)]$A^G$ contains a separable element of
    degree $n$ over $\Z(A)$.
  \item[(c)]For $x\in X$ in general position, $\Stab_{G \times
      \PGLn}(x)$ is contained in $G \times T_x$, where $T_x$ is
      a torus in $\PGLn$.
\qed
  \end{itemize}
\end{cor}

Here by a separable element of $A$ we mean an element whose minimal
polynomial over $\Z(A)$ is separable, i.e., has no repeated roots.

\begin{remark} \label{rem.m2}
  It is necessary in Corollary~\ref{cor.of.thm2}(b) to require that
  $a$ is separable over $\Z(A)$.  Indeed, in Example~\ref{ex:m2} set
  $n = 2$, $A = \Mat_2(k)$ and $G =
  \{\bigl(\begin{smallmatrix}1&g\\0&1\end{smallmatrix}\bigr) \, | \, g
  \in k \} $.  Then the fixed algebra $A^G$ consists of all matrices
  of the form
  $\bigl(\begin{smallmatrix}a&b\\0&a\end{smallmatrix}\bigr)$ with
  $a,b\in k$. In particular, $A^G$ contains elements of degree $n=2$
  over $\Z(A)=k$, but the minimal polynomial of any such element has
  repeated roots.
\end{remark}

\section{The $G$-action on the center of $A$}
\label{sect.center}

Throughout this section, we consider a geometric action of an
algebraic group $G$ on a central simple algebra~$A$ of degree~$n$ with
associated $G\times\PGLn$-variety $X$.  It is sometimes possible to
deduce information about the $G$-action on $A$ from properties of the
$G$-action on the center $\Z(A)$.  In this section, we find conditions
on the $G$-action on $\Z(A)$ which allow us to answer question (a)
in~\eqref{q1}.  

Recall that the field of rational functions on $X/\PGLn$ is
$G$-equivariantly isomorphic to the center $\Z(A)$ of $A$ (see
Lemma~\ref{lem4.1-NEW}).  Of course, a priori $X/\PGL_n$ is only
defined up to birational isomorphism. From now on we will fix a
particular model $W$ equipped with a regular $G$-action and a
$G$-equivariant rational quotient map for the $\PGLn$-action on $X$
\[
  \pi \colon X \dasharrow W\,.
\]
It will not matter in the sequel which model $W$ of $X/\PGLn$ we use.
Note that the $G$-variety $W$ is just a birational model
for the $G$-action on $\Z(A)$.  In many (perhaps, most) cases,
$W$ is much easier to construct than $X$; for an example
of this phenomenon, see Section~\ref{sect.ex-mult}.

We begin with a simple observation, relating stabilizers in $X$ and $W$.

\begin{lem} \label{lem7.1}
  Let $X$ be a $G\times\PGLn$-variety which is $\PGLn$-generically
  free.  Denote by $\pi\colon X\dasharrow X/\PGLn$ the rational
  quotient map for the $\PGLn$-action.  Then for $x\in X$ in general
  position, the projection $G \times \PGLn \lra G$ onto the
  first factor induces an isomorphism between $\Stab_{G \times
    \PGLn}(x)$ and $\Stab_G(\pi(x))$.
\end{lem}

\begin{proof}
  For $x \in X$ in general position, $\pi$ is defined at $x$, the
  fiber over $\pi(x)$ is the orbit $\PGLn x$, and $\Stab_{\PGLn}(x)$
  is trivial. For such $x$, the projection~$p$ restricts to a
  surjective map
  \[ \Stab_{G \times \PGLn}(x) \lra \Stab_G(\pi(x)) \]
  whose kernel is $\Stab_{\PGLn}(x) = \{ 1 \}$, and the lemma follows.
\end{proof}

\begin{prop} \label{prop7.2}
  \begin{itemize}
  \item[\textup{(a)}] Suppose that for $w \in W$ in general position,
    the stabilizer $\Stab_G(w)$ does not admit a non-trivial
    homomorphism to $\PGLn$.
    Then $A^G$ is a central simple algebra of degree $n=\deg (A)$.

\smallskip
  \item[\textup{(b)}] Suppose that for $w \in W$ in general position,
   $\Stab_G(w)$ is an abelian group consisting of semisimple elements and
   the $n$-torsion subgroup of $\Stab_G(w)/\Stab_G(w)^0$ is cyclic.
   Then there exists an $a \in A^G$ with $n$ distinct eigenvalues.
  \end{itemize}
\end{prop}

Note that the condition of part (a) is satisfied if the $G$-action on
$W$ is generically free.

\begin{proof}
  (a) By Lemma~\ref{lem7.1}, $\Stab_{G \times
    \PGLn}(x)\subseteq G\times\{1\}$ for $x$ in general
  position in $X$.  The desired conclusion follows from
  Theorem~\ref{thm2}(a).

  \smallskip (b) Let $H$ be the projection of $\Stab_{G \times
    \PGLn}(x)$ to $\PGLn$.  By Lemma~\ref{lem7.1}, $H$ is an abelian
  group consisting of semisimple elements, and $H/H^0$ is a homomorphic
  image of $\Stab_G(w)/\Stab_G(w)^0$.  Using the fundamental theorem
  of finite abelian groups, one checks that surjective homomorphisms
  of finite abelian groups preserve the property that the $n$-torsion
  subgroup is cyclic.  By~\cite[Corollary 2.25(a)]{steinberg}, $H$ is
  contained in a maximal torus of $\PGLn$.  (Note that the torsion
  primes for $\PGLn$ are the primes dividing~$n$;
  see~\cite[Corollaries 1.13 and 2.7]{steinberg}.)  The desired
  conclusion now follows from Theorem~\ref{thm2}(b).
  \end{proof}

We will now use Proposition~\ref{prop7.2}, to study inner actions.
Recall that an automorphism $\phi$ of a central simple
algebra $A$ is called {\em inner} if there exists an invertible
element $a \in A$ such that $\phi(x) = ax a^{-1}$ for every $x \in A$,
and {\em outer} otherwise. By the Skolem-Noether theorem $\phi$
is inner if and only $\phi(x) = x$ for every $x \in \Z(A)$.

\begin{cor} \label{cor4}
  Let $G$ be a finite group or a torus acting geometrically
  on a central simple
  algebra $A$ of degree $n$. The elements of $G$ that act by inner
  automorphisms form a normal subgroup of $G$; denote this subgroup by
  $N$.
  \begin{itemize}
  \item[\textup(a)]If $N = \{ 1 \}$ \textup{(}i.e., if $G$ acts on $A$ by
    outer automorphisms\textup{)}, then $A^G$ is a central simple
    algebra of degree $n$.
  \item[\textup(b)]If $N$ is a cyclic group, then there is an element
  $a \in A^G$ with $n$ distinct eigenvalues.
  \end{itemize}
\end{cor}

 In the case where the group $G$ is finite, part (a) is proved by
 algebraic means and under weaker hypotheses in~\cite[Theorem~2.7 and
 Corollary 2.10]{montgomery}.  Note also that since
 every action of a finite group on a central simple algebra 
 is algebraic (see Definition~\ref{defn:algebraic-action}), our 
 assumption that the action is geometric is only relevant if $G$ is a torus.
 Moreover, if $G$ is a torus then every geometric action is algebraic;
 see Corollary~\ref{cor:prop:geom=>alg}.
 
\begin{proof} We may assume that the action is faithful. Indeed,
if $K \subseteq G$ is the kernel of this action, we can replace
$G$ by $G/K$ and $N$ by $N/K$.

Now let $W$ be an irreducible $G$-variety whose function field $k(W)$
is $G$-equivariantly isomorphic to $\Z(A)$ (over $k$); see the beginning
of Section~\ref{sect.center}.
Clearly an element of $G$ acts trivially on $W$ if
and only if it acts on $A$ by an inner automorphism. Now recall that
if $G$ is a finite group or a torus then the stabilizer in general
position for the $G$-action on $W$ is precisely the kernel~$N$ of this
action; cf.\ Lemma~\ref{lem.diag}.

The desired conclusions in parts (a) and (b) now follow from parts (a)
and (b) of Proposition~\ref{prop7.2}, respectively.
\end{proof}

\section{Which geometric actions are algebraic?}
\label{sect.geom=>alg}

Theorem~\ref{thm1}(a) says that every algebraic action is
geometric. It is easy to see that the converse is not true.  For
example, let $Y$ be a generically free $\PGLn$-variety (e.g.,
we can take $Y = \PGLn$ where $\PGLn$ acts on $Y$ by translations),
and consider the $G \times \PGLn$-variety $X = (G/P) \times Y$,
where $G$ is a non-solvable connected algebraic group, and $P$
is a proper parabolic subgroup. Here $G$ acts by translations on the
first factor, and $\PGLn$ acts on the second factor.  
Since the $\PGLn$-action on $X$ is
generically free, $A = k_n(X)$ is a central simple
algebra of degree~$n$.  On the other hand, since $G/P$ is
complete, it is easy to see that $X$ is not birationally isomorphic to
an affine $G \times \PGLn$-variety; hence by
Theorem~\ref{thm1}(b), this action is not algebraic. 

Nevertheless, we will now show that under fairly mild
assumptions, the converse of Theorem~\ref{thm1}(a) holds,
i.e., every geometric action is, indeed, algebraic.

\begin{lem} \label{lem:geom=>alg}
  Let $G$ be an algebraic group, and let $X$ be an irreducible
  $G\times\PGLn$-variety which is $\PGLn$-generically free.
  Assume that $X$ has a stable affine model as a $G \times \PGLn$-variety.
  Then the induced action of $G$ on $k_n(X)$ is     algebraic.
\end{lem}

\begin{proof}
  We may assume without loss of generality that $X$ itself is affine
  and stable as a $G \times \PGL_n$-variety. By
  Theorem~\ref{thm1}(b) it suffices to show that $X$ is stable 
  as a $\PGLn$-variety, i.e., that $\PGLn$-orbits in general
  position in $X$ are closed. Let $x \in X$ be a point in general position.
  Then the $G \times \PGLn$-orbit $(G \times \PGLn) x$ is closed 
  in $X$ and can be naturally identified with the homogeneous space 
  $(G \times \PGLn)/H$, where $H =  \Stab_{G \times \PGLn}(x)$.
  The $\PGLn$-orbit $(\PGLn) \cdot x$ is then
  identified with the image $Z$ of $\PGLn$ in $(G \times \PGLn)/H$. 
  It thus remains to show that $Z$ is closed in $(G\times\PGLn)/H$.
  Indeed, $Z$ is also the image of the product $\PGLn H$, which 
  is a closed subgroup of $G\times\PGLn$ (because $\PGLn$ is normal; see
  \cite[\S7.4]{humphreys}).  Since $\PGLn H$ is a closed subgroup of
  $G \times \PGLn$ containing $H$, its image $Z$ in $(G\times\PGLn)/H$ 
  is closed; see \cite[\S12.1]{humphreys}.
\end{proof}

 \begin{cor} \label{cor2:geom=>alg}
   Let $G$ be an algebraic group, and let $X$ be an irreducible
   $G\times\PGLn$-variety which is $G\times\PGLn$-generically free.
   Then the induced action of $G$ on $k_n(X)$ is algebraic.
 \end{cor}

 \begin{proof}
 By~\cite[Theorem~1.2(i)]{rv3} $X$ has a stable affine birational model as
 a $G\times\PGLn$-variety. The desired conclusion is now
 immediate from Lemma~\ref{lem:geom=>alg}.
 \end{proof}

The criterion for a geometric action to be algebraic
given by Lemma~\ref{lem:geom=>alg} can be further simplified
by considering the $G$-action on the center of $A$, as
in Section~\ref{sect.center}.

\begin{prop} \label{prop:geom=>alg}
  Consider a geometric action of an algebraic group $G$ on a
  central simple algebra $A$, and let $W$ be a birational model for
  the $G$-action on $\Z(A)$.  Then the $G$-action on $A$ is
  algebraic, provided one of the following conditions holds:
  \begin{itemize}
  \item[\textup{(a)}]The $G$-action on $W$ is generically free.
  \item[\textup{(b)}]The normalizer $H = \N_G(G_w)$ is reductive for $w
    \in W$ in general position.
  \item[\textup{(c)}]$G$ is reductive and the stabilizer $G_w$ is
    reductive for $w \in W$ in general position.
  \item[\textup{(d)}]$G$ is reductive and $W$ has a stable affine
    model as $G$-variety.
  \end{itemize}
\end{prop}

\begin{proof} Let $X$ be an associated $G \times \PGLn$-variety
  for the $G$-action on $A$. Recall that the $\PGLn$-action on $X$ is
  generically free and $W$ is the rational quotient $X/\PGLn$.  In
  view of Lemma~\ref{lem:geom=>alg}, it suffices to show that $X$
  has a stable affine model as a $G \times \PGLn$-variety.

  \smallskip (a) Immediate from Corollary~\ref{cor2:geom=>alg} and
  Lemma~\ref{lem7.1}.

  \smallskip (b) Choose $x \in X$ in general position, and set $w =
  \pi(x) \in W$.  Let $S_x = \Stab_{G \times \PGLn}(x)$.  We claim that
  $\N_{G \times \PGLn}(S_x)$ is reductive for $x \in X$ in general
  position. The desired conclusion follows from this claim
  by~\cite[Theorem 1.2(ii)]{rv3}.
  
  The proof of the claim is based on two simple observations.  First
  of all, if $H = \N_G(G_w)$ is reductive, then so is $S_x\simeq G_w$.
  Indeed, the unipotent radical of $\operatorname{R}_u(G_w)$ is
  characteristic in $G_w$, hence, normal in $H$. Since $H$ is
  reductive, this implies $\operatorname{R}_u(G_w) = \{ 1 \}$, i.e.,
  $G_w$ is reductive, as claimed.
  
  Secondly, by Lemma~\ref{lem7.1}, the normalizer $\N_{G \times
    \PGLn}(S_x)$ is a priori contained in $H \times \PGLn$, i.e.,
  \[ \N_{G \times \PGLn}(S_x) = \N_{H \times \PGLn}(S_x) \, . \]
  Since both $H \times \PGL_n$ and $S_x$ are reductive, the normalizer
  $\N_{H \times \PGLn}(S_x)$ is reductive as well; see~\cite[Lemma
  1.1]{lr}.  This concludes the proof of the claim and thus of part
  (b).

  \smallskip (c) If $G$ and $G_w$ are both reductive then using
  \cite[Lemma 1.1]{lr} once again we see that $\N_G(G_w)$ is also
  reductive. Part~(c) now follows from part~(b).

  \smallskip (d) After replacing $W$ by a stable affine model, we see
  that for $w\in W$ in general position, the orbit $Gw\simeq G/G_w$ is
  affine, so that $G_w$ is reductive by Matsushima's theorem,
  see~\cite[Theorem~4.17]{pv}.  Now use part~(c).
\end{proof}

\begin{cor} \label{cor:prop:geom=>alg}
  Let $G$ be an algebraic group whose connected component is a
  torus.  Then every geometric action of $G$ on a central simple
  algebra is algebraic.
\end{cor}

\begin{proof}
  In this case, every subgroup of $G$ is reductive, so that part~(c)
  of Proposition~\ref{prop:geom=>alg} applies.
\end{proof}

\section{Proof of Theorem~\ref{thm4}} 
\label{sect.thm4}

\numberedpar {\em The generic torus.}\label{generic.torus}
Let $T$ be a maximal torus in $\GL_n$, and let $N$ be the normalizer of
the image of $T$ in $\PGLn$.  Since $\PGLn$ permutes the maximal tori
in $\GL_n$ transitively, one can think of $\PGLn/N$ as the variety of
maximal tori of $\GL_n$ (or equivalently, of $\PGLn$).  We briefly
recall how one can construct a $\PGLn$-equivariant rational map
\[ \pi \colon \Mn \dasharrow \PGLn/N \]
which sends a non-singular matrix $\alpha \in \Mn$ with distinct
eigenvalues to the unique maximal torus in $\GL_n$ containing
$\alpha$.  The map $\pi$ is sometimes called the {\em generic torus}
of $\GL_n$; cf.~\cite[4.1]{voskresenskii}.

Denote by $\Gr(n,n^2)$ the Grassmannian of $n$-dimensional subspaces
of $\Mn$.  The action of $\PGLn$ on $\Mn$ induces a regular action of
$\PGLn$ on $\Gr(n,n^2)$.  Define a rational, $\PGLn$-equivariant map
$\pi_1 \colon\Mn\dasharrow \Gr(n,n^2)$ by sending a non-singular matrix
$\alpha$ with distinct eigenvalues to
$\Span(1,\alpha,\ldots,\alpha^{n-1})$.  The unique maximal torus
$T_{(\alpha)}$ of $\GL_n$ containing $\alpha$ is characterized by
$\Span(T_{(\alpha)})=\pi_1(\alpha)$.  The image of $\pi_1$ consists
thus of a single $\PGLn$-orbit~$O$.  
Since the stabilizer of both $T$
and $\Span(T)$ is $N$, $gN\mapsto g\Span(T)g\inv$ defines an
isomorphism $\pi_2\colon \PGLn/N\to O$. 
Here $T$ is the maximal torus in $\GL_n$ which we chose (and fixed)
at the beginning of this section and $N$ is the normalizer 
of the image of $T$ in $\PGLn$. Now 
$\pi=\pi_2\inv\circ\pi_1$ is a $\PGLn$-equivariant rational map $\Mn
\dasharrow \PGLn/N$ 
such that for any $\alpha$ as above, $\pi(\alpha)=gN$ if and only if
$gTg\inv$ is the unique torus of $\GL_n$ containing $\alpha$.

\smallskip

\numberedpar {\em Proof of Theorem~\ref{thm4}.}
(a) Suppose $A=k_n(X)$ has a $G$-invariant maximal
\'etale subalgebra $E$.  It follows easily from the primitive element
theorem that there is an $a \in E$ so that $E = \Z(A)[a]$.  Choose 
one such $a$.  By
Lemma~\ref{lem:etale}, $a$ has distinct eigenvalues.  Adding some
constant in $k$ to $a$, we may assume that the eigenvalues of $a$ are
distinct and nonzero.  Then for $x \in X$ in general position $a(x)$
is a matrix whose eigenvalues are distinct and nonzero.  We now define
a rational map $\varphi\colon X \dasharrow \PGLn/N$ by
$\varphi(x)=\pi(a(x))$.  This map is $\PGLn$-equivariant by
construction.  Moreover, for every $g \in G$, $g(a) \in E$ commutes
with $a$.  Thus, for $x \in X$ in general position, $a(x)$ and
$g\inv(a)(x)=a(g(x))$ lie in the same maximal torus, and consequently,
$\varphi(x)=\varphi(g(x))$.

Conversely, suppose there exists a $G \times \PGLn$-equivariant
rational map $X \dasharrow \PGLn/N$. After removing the indeterminacy
locus from $X$, we may assume this map is regular. We may
also assume that $\PGLn$ acts freely on $X$.  Let $X_0$ be the
preimage of the coset $N\in\PGLn/N$ in $X$. Note that $X_0$ is $G
\times N$-invariant, that $X=\PGLn\cdot X_0$, and that the $N$-action
on $X_0$ is generically free. Moreover, $X$ is birationally isomorphic
as $\PGLn$-variety to $\PGLn * _N X_0$, see~\cite[Theorem 1.7.5]{popov2}.

Let $\Delta \simeq \bbA^n$ be the variety of diagonal $n \times n$-matrices.
By~\cite[Proposition 7.1]{reichstein} there exists an $N$-equivariant
rational map $a \colon X_0 \dasharrow \Delta$ whose image
contains a matrix
with distinct eigenvalues. (Note that here we use the fact that
$\Delta$ is a vector space and $N$ acts on it linearly.)
This rational map then naturally extends to
a $\PGLn$-equivariant rational map
\[ X \simeq \PGLn *_N X_0 \dasharrow \PGLn *_N \Delta \simeq \Mn  \]
induced by $(g, x_0) \mapsto (g, a(x_0))$.
By abuse of notation, we denote this extended rational
map by $a$ as well.

We now view $a$ as an element of $A=k_n(X)$.
Since the image of $a$ contains a matrix with distinct eigenvalues,
Lemma~\ref{lem:etale} tells us that $E = \Z(A)[a]$ is a maximal \'etale
subalgebra of $A$. It remains to show that $E$ is $G$-invariant.
To do this it suffices to prove that $g(a) \in E$ for every $g \in G$.
Since $E = C_A(E)$, we only need to establish that $g(a)$ commutes
with $a$, i.e., that the commutator
$b = [a, g(a)]$ equals $0$. Indeed, for any $x \in X_0$,
\[ b(x) = [a(x), a(g^{-1}(x))] = [a(x), a(y)] \, , \]
where $y = g^{-1}(x) \in X_0$. By our construction $a$ maps every
element of $X_0$ to a diagonal matrix. In particular, $a(x)$ and
$a(y)$ commute, and thus $b(x) = 0$ for every $x \in X_0$. Since $b$
is a $\PGLn$-equivariant rational map $X \dasharrow \Mn$ and since
$\PGLn \cdot X_0 =X$, we conclude that $b = [a, g(a)]$ is identically
zero on $X$, as claimed. This completes the proof of part (a).

\smallskip
(b) The action of $G \times \PGLn$ on $\PGLn/N$ has stabilizer
of the form $G \times \N(S)$ at every point, where $S$ is a maximal
torus of $\PGLn$.  Part (b) is now an immediate consequence
of part (a).

\smallskip (c) Assume that $A$ has a $G$-invariant maximal \'etale
subalgebra.  Let $x \in X$ be a point in general position. We claim
that
\begin{equation} \label{e7.1}
\dim(Gx \cap \PGLn x) \leq n -1 \, .
\end{equation}
Indeed, $Gx \cap \PGLn x$ is easily seen to be the image of the
morphism from $\Stab_{G \times \PGLn}(x)$ to $X$ given by $(g, p)
\mapsto px$.  Since $\Stab_{G \times \PGLn}(x) \subseteq G \times
\N(T_x)$ by part (b), we conclude that
\[ \dim(Gx \cap \PGLn x) \leq \dim \, \N(T_x) = n-1 \, \]
as claimed.

Consider the rational quotient map $\pi \colon X \dasharrow
X/\PGLn$.  We may assume without loss of generality that $\pi$ is
defined at $x$.  Now restrict $\pi$ to the (well-defined) rational map
$\pi_{G^0x} \colon G^0x \dasharrow X/\PGLn$, where $G^0$ is
the connected component of $G$.
For $y\in G^0x$ in general position, the fiber over $\pi_{G^0x}(y)$ is
$G^0x\cap\PGLn y=G^0y\cap\PGLn y$.
By~\eqref{e7.1},
\begin{align*}
\dim(Gx) = \dim(G^0x) &\le \dim(X/\PGLn) + n - 1 \\
         &= \dim(X) - \dim(\PGLn) + n - 1\\ & = \dim(X) - n^2 + n \, .
\end{align*}
So $\dim(X)-\dim(Gx)\geq n^2-n$.  This proves part (c).
\qed

\section{Proof of Theorem~\ref{thm3}}
\label{sect.thm3}

We begin by spelling out what it means for an algebraic group 
action on a central simple algebra to be split 
in terms of the associated variety.

\begin{lem} \label{lem.G-split}
A geometric action of an algebraic group $G$ on a central simple algebra
$A$ of degree $n$ is $G$-split in the sense of Definition~\ref{def.G-split} 
if and only if its associated $G \times \PGLn$-variety is birationally
isomorphic to $X_0 \times \PGLn$, for some $G$-variety $X_0$.
\end{lem}

Here $G$ acts on the first factor and $\PGLn$ acts on the second
factor by translations.

\begin{proof}
Suppose $X = X_0 \times \PGLn$. Then
we have the following $G$-equivariant isomorphisms, 
\[ \RMaps_{\PGLn}(X, \Mn) \simeq \RMaps(X_0, \Mn) \simeq 
\Mn(k) \otimes_k k(X_0) \, , \]
where the first isomorphism is given by $f \mapsto f|_{X_0 \times 1_{\PGLn}}$
for every $\PGLn$-equivariant rational map $f \colon X \dasharrow \Mn$.
In other words, the induced $G$-action on $A = k_n(X)$ is $G$-split 
in the sense of Definition~\ref{def.G-split}.

Conversely, suppose a geometric $G$-action on $A$ is $G$-split.
Denote the associated $G \times \PGLn$-variety by $X$. Let $X_0 = X/\PGLn$
be the rational quotient of $X$ by the $\PGLn$-action. Note that
$k(X_0) = \Z(A)$. Then, as we saw above,
$\RMaps_{\PGLn}(X_0 \times \PGLn,\Mn)$ is $G$-equivariantly isomorphic to
$\Mn(k) \otimes_k k(X_0)$, which is $G$-equivariantly 
isomorphic to $A$ (because $A$ is $G$-split).  
By Corollary~\ref{cor:unique}, we conclude 
that $X$ is birationally isomorphic to $X_0 \times \PGLn$.
\end{proof}

\begin{cor} \label{cor.G-split}
  Consider a geometric action of an algebraic group $G$ on a central
  simple algebra $A$, with associated $G \times \PGLn$-variety $X$.
  Then for any $G$-variety $X_0$ the following are equivalent:
  \begin{itemize}
  \item[\upshape(a)]$L = k(X_0)$ is a $G$-splitting field for $A$.
  \item[\upshape(b)]There exists a dominant rational map $f\colon X_0
    \times \PGLn \dasharrow X$ which is $G\times\PGLn$-equivariant.
\end{itemize}
\end{cor}

Here $G$ acts on the first factor of $X_0 \times \PGLn$ and $\PGLn$
acts on the second factor by translations, as in
Lemma~\ref{lem.G-split}.

\begin{proof} (a) $\Longrightarrow$ (b). The $G$-action on
  $A' = A \otimes_{\Z(A)} L \simeq \Mn(k) \otimes_k L$ is geometric,
  with associated variety $X' = X_0 \times \PGLn$; see
  Lemma~\ref{lem.G-split}.  The embedding $j \colon A \hookrightarrow
  A'$ induces a $G\times\PGLn$-equivariant dominant rational map $j_*
  \colon X' \dasharrow X$; see Lemma~\ref{lem:unique}.
  
  \smallskip (b) $\Longrightarrow$ (a): Let $X' = X_0 \times \PGLn$.
  By Lemma~\ref{lem:unique}, $f$ induces a $G$-equivariant embedding
  $f^* \colon A \hookrightarrow A'$ of central simple algebras, where
  $A' = k_n(X') \simeq \Mn (k) \otimes_k k(X_0)$; see
  Lemma~\ref{lem.G-split}. In other words, $A'$ is $G$-equivariantly
  isomorphic to $A \otimes_{\Z(A)} k(X_0)$.
\end{proof}

\numberedpar {\em Proof of Theorem~\ref{thm3}}.
Let $X$ be the associated $G\times\PGLn$-variety for the $G$-action on
$A$.  Consider the dominant morphism
$f \colon X \times \PGLn \lra X$ given by $(x, h) \mapsto h x$.
If we let $(g, h) \in G \times \PGLn$ act on $X \times \PGLn$ 
by $(g, h) \cdot (x, h') = (g x, h h')$, as in Lemma~\ref{lem.G-split}
and Corollary~\ref{cor.G-split}, then we can easily check that 
$f$ is $G \times \PGLn$-equivariant. 
By Corollary~\ref{cor.G-split}, we conclude that
$L = k(X)$ is a $G$-splitting field for $A$.
Moreover, 
\begin{align*} 
    \trdeg_{\Z(A)}  L &= \trdeg_k (L) - \trdeg_k \Z(A)\\
                      &=\dim(X) - \dim(X/\PGLn) = n^2 - 1,
  \end{align*}
as claimed.  Note that if $G$ acts algebraically on $A$, we may assume
that $X$ is affine by Theorem~\ref{thm1}(b).
\qed

\section{More on $G$-splitting fields}
\label{sect3b}

In this section we discuss $G$-splitting fields 
in the case where $G$ is a connected group. Our main result 
is the following:

\begin{prop} \label{prop3} Consider a geometric action of
  a connected algebraic group $G$ on a central simple algebra $A$ of
  degree $n$.
  Then there exists an affine $G$-variety $X_0$ such that $L = k(X_0)$
  is a $G$-splitting field of $A$ and
  \begin{equation} \label{e.statement.prop3} 
    \trdeg_{\Z(A)} \, L = \dim \, \Stab_{G \times \PGLn}(x)
                        = \dim \, \Stab_G (w)   \, , 
  \end{equation}  
  where $x$ and $w$ are points in general position in the associated
  $G \times \PGLn$-variety~$X$ and in the rational quotient
  $W=X/\PGLn$, respectively.  In particular,
  \[\trdeg_{\Z(A)} \, L \leq \dim(G)\,.\]
\end{prop}

\medskip

Note that for $w \in W$ in general position we have
\begin{equation}\label{e.prop3}
\begin{split}
\dim \, \Stab_G(w)  & = \dim(G) - \dim(Gw)   \\ 
                   & = \dim(G) - \bigl(\dim(W) - \dim(W/G)\bigr) \,.
\end{split}
\end{equation}
so that the integer  $\dim \, \Stab_G (w)$ for $w \in W$ in
general position, which appears in the statement 
of Proposition~\ref{prop3} is well defined. 
Similarly, the integer
$\dim \, \Stab_{G \times \PGLn}(x)$ for $x \in X$ in general 
position is also well-defined.
Since $\trdeg_{\Z(A)^G} \, \Z(A)=\dim(W)-\dim(W/G)$,
\eqref{e.statement.prop3} can be restated in algebraic terms as
\addtocounter{equation}{-1}
\begin{equation}\renewcommand{\theequation}{\ref{e.statement.prop3}$'$}
      \trdeg_{\Z(A)}  L  =  \dim(G) - \trdeg_{\Z(A)^G} \, \Z(A) \, . 
\end{equation}

In general, the value for $\trdeg_{\Z(A)} L$ given in
\eqref{e.statement.prop3} and (\ref{e.statement.prop3}$'$) is the
smallest possible, see Remark~\ref{rem.prop3.minimal}.  
Our proof of Proposition~\ref{prop3} will rely on the following lemma.

\begin{lem} \label{lem.cover} Let $H$ be a connected algebraic group and
  let $V$ be an irreducible $H$-variety.  Then there
  exists an irreducible variety $Y$ and an $H$-equivariant dominant
  morphism $Y \times H \lra V$ such that $\dim(Y) = \dim(V/H)$.
\end{lem}

The action of $H$ on $Y\times H$ is induced by the trivial action on
$Y$ and by the translation action on $H$.

\begin{proof} See~\cite[(1.2.2)]{popov2} or~\cite[Proposition 2.7]{pv},
where the term quasi-section is used to describe $Y$.
\end{proof} 

\numberedpar {\em Proof of Proposition~\ref{prop3}}.
By Lemma~\ref{lem.cover} (with $H = G \times \PGLn$)
there is a $G \times \PGLn$-equivariant dominant morphism 
$f \colon Y \times (G \times \PGLn) \lra X$, where
\begin{equation} \label{e.dim1}
\dim(Y) = \dim\bigl(X/(G \times \PGLn)\bigr) = \dim(W/G)\, .
\end{equation}
Note that since $G \times \PGLn$ acts trivially 
on $Y$, we can take $Y$ to be affine.  Setting 
$X_0 = Y \times G$ (as a $G$-variety) and applying 
Corollary~\ref{cor.G-split}, we conclude 
that $L = k(X_0)$ is a $G$-splitting field for $A$.
By our construction, $X_0 = Y \times G$ is affine.
Since the second equality in~\eqref{e.statement.prop3}
is an immediate consequence of Lemma~\ref{lem7.1},
we only need to check that $\trdeg_{\Z(A)} \, L = \dim \, \Stab_G(w)$
for $w \in W$ in general position.  Indeed,
\begin{align*} 
    \trdeg_{\Z(A)} \, L &= \trdeg_k(L) - \trdeg_k \Z(A)
                     =\dim (X_0) - \dim(X/\PGLn) \\
                     & = \dim(Y) + \dim(G) - \dim(X/\PGLn) \\  
                     & = \dim(G) - \bigl(\dim(W) - \dim(W/G)\bigr)
                       = \dim \, \Stab_G(w) \, ,
\end{align*}
where the two last equalities follow from~\eqref{e.dim1}
and~\eqref{e.prop3}, respectively.
\qed\medskip

Specializing Proposition~\ref{prop3} to the case of torus actions, we
recover a result which was proved in \cite{vonessen:solv} for
algebraic actions in arbitrary characteristic.

\begin{cor} \label{cor3}
Suppose a torus $T$ acts geometrically {\upshape(}or equivalently, 
algebraically; cf.\ Corollary~\ref{cor:prop:geom=>alg}{\upshape)} 
on a central simple algebra $A$.  Let $H$ be the kernel 
of the $T$-action on $\Z(A)$.  Then there exists 
a $T$-variety $X_0$ such that $L = k(X_0)$
is a $T$-splitting field for $A$ and $\trdeg_{\Z(A)} \, L = \dim(H)$.
\end{cor}

\begin{proof} Let $X$ be the associated $G \times \PGLn$-variety
and $W = X/\PGLn$, as before. 
  By Lemma~\ref{lem.diag}, applied to the $T$-action on
  $W$, we have $H = \Stab_T(w)$ for $w \in W$ in general position.  
  The corollary now follows from Proposition~\ref{prop3}.
\end{proof}
 
\begin{remark}\label{rem.prop3.minimal}
  If the $T$-action on $A$ is faithful
  then the value of $\trdeg_{\Z(A)} \, L$ given by Corollary~\ref{cor3}
  is the smallest possible. Indeed, since the $T$-action on both $A$
 and $L = k(X_0)$ is algebraic (cf.\  Corollary~\ref{cor:prop:geom=>alg}), 
  \cite[Theorem~2(b)]{vonessen:solv} tells us 
  that $\trdeg_{\Z(A)} \, L \ge \dim(H)$ 
  for every $T$-splitting field of the form $L = k(X_0)$,
  where $X_0$ is a $T$-variety. 
  \qed
\end{remark}

\begin{remark} Suppose a torus $T$ acts geometrically 
(or equivalently, algebraically; cf.\ Corollary~\ref{cor:prop:geom=>alg}) 
on a division algebra $D$.  Then \cite[Theorem 2(c)]{vonessen:solv} 
asserts that $D$ has a $T$-splitting field $L$ of the form $k(X_0)$ 
such that $[L: \Z(D)] < \infty$.

We now give an alternative proof of this result based on
Corollary~\ref{cor3}. Let $T_0 \subset T$ be the kernel 
of the $T$-action on $D$. After replacing $T$ by $T/T_0$, 
we may assume the action is faithful. 
Let $H$ be the kernel of the $T$-action on $\Z(A)$,
i.e., the subgroup of $T$ acting by inner automorphisms.
By Corollary~\ref{cor.inner.3}, $H$ is a finite group.
By Corollary~\ref{cor3}, there exists a $T$-splitting field $L =
k(X_0)$ such that $\trdeg_{\Z(D)} \, L = \dim(H) = 0$.
Since $L$ is finitely generated over $k$ (and hence, over $\Z(D)$), 
we conclude that $[L : \Z(D)] < \infty$.
\qed
\end{remark}

\section{An example: algebraic actions of unipotent groups} 
\label{sect.unipotent}

In this and the subsequent three sections we will present
examples, illustrating Theorems~\ref{thm2}, \ref{thm4}, and~\ref{thm3}.
We begin by applying Theorems~\ref{thm2} and~\ref{thm4} in
the context of unipotent group actions on division algebras.

\begin{prop} \label{prop.unipotent} 
  Let $U$ be a unipotent group acting algebraically
  on a finite-dimensional division algebra $D$. Then
  $D^U$ is a division algebra of the same degree as $D$.
\end{prop}

\begin{proof}
  Say $D$ has degree~$n$, and let $X$ be the associated
  $U\times\PGLn$-variety.  By Lemma~\ref{lem7.1}, for $x \in X$ in
  general position, $\Stab_{U \times \PGLn}(x)$ is a unipotent group
  (it is isomorphic to a subgroup of $U$).  Consequently, the
  projection $H_x$ of this group to $\PGLn$ is unipotent.
  
  On the other hand, by~\cite[Proposition 7]{vonessen:solv}, $D$ has a
  $U$-invariant maximal subfield. In view of Theorem~\ref{thm4}(b),
  this implies that $H_x$ is a subgroup of the normalizer of a maximal
  torus in $\PGLn$; in particular, $H_x$ has no non-trivial unipotent
  elements.  This is only possible if $H_x = \{ 1 \}$, i.e., if
  \[ \Stab_{U \times \PGLn}(x) \subseteq U \times \{ 1 \} \, . \]
  The desired conclusion now follows from Theorem~\ref{thm2}(a).
\end{proof}

\begin{remark} The condition that $D$ is a division algebra 
  is essential here. Suppose $G = U$ is a non-trivial unipotent
  subgroup of $\PGLn$ acting on $A = \Mn(k)$ by conjugation, as in
  Example~\ref{ex.conjugation}. Since $A$ is a finite-dimensional
  $k$-vector space, this action is easily seen to be algebraic. On the
  other hand, the fixed algebra $A^U$ is not a central simple algebra
  of degree $n$; cf.\ Example~\ref{ex:m2}(a) (see also
  Remark~\ref{rem.m2}).
\end{remark} 

\section{An example: the $\GL_m$-action on $\UDmn$, $m \ge n^2$}
\label{sect.ex.UDmn}

We now return to the $\GL_m$-action on the universal division 
algebra $A = \UDmn$, described in Example~\ref{ex.universal}. 
In this section we will assume that $m \ge n^2$; in the next 
section we will set $m = n = 2$. The case where $m \le n^2 - 1$ will be
considered in~\cite{rv5}. 

\begin{prop} \label{prop.split} Let  $A = \UDmn$, where $m \ge n^2$.
Then
\begin{itemize}
\item[\upshape(a)] $A^{\GL_m} = k$
\item[\upshape(b)] $\trdeg_{\Z(A)} \, L \ge n^2 - 1$ for every
  $\GL_m$-splitting field $L$ of $A$ of the form $L = k(X_0)$, where
  $X_0$ is a $\GL_m$-variety.
\end{itemize}
\end{prop}

Part (b) shows that the value of $\trdeg_{\Z(A)} \, L$
given by Theorem~\ref{thm3} is optimal for this action.

\begin{proof}
  The variety $\Mnm$ is an associated $\GL_m
  \times \PGLn$-variety for the $\GL_m$-action on~$A$; see
  Example~\ref{ex.universal}. The key fact underlying the proof of
  both parts is that for $m\geq n^2$, $\Mnm$ has 
  a dense $\GLm$-orbit; denote this orbit by~$X$.
  Since the actions of $\GLm$ and $\PGLn$ commute,
  $X$ is $\PGLn$-stable, and therefore is also 
  an associated $\GL_m \times \PGLn$-variety for the $\GL_m$-action on~$A$.

\smallskip
(a) By Remark~\ref{rem:first-results}, 
$A^{\GL_m} = \RMaps_{\PGLn}(X/\GL_m, \Mn)$. Since
$X$ is a single $\GL_m$-orbit, the rational quotient $X/\GL_m$ is
a point (with trivial $\PGLn$-action).  Clearly, every 
$\PGLn$-equivariant rational map $f \colon \{ pt \} \dasharrow \Mn$
is regular and has its image in the center of $\Mn$. In other words,
\[ A^{\GL_m} = \RMaps_{\PGLn}(X/\GL_m, \Mn) = \RMaps(\{ pt \}, k) = k \, , \]
as claimed.

\smallskip
(b) By Corollary~\ref{cor.G-split} there exists a dominant rational map
  $f \colon X_0 \times \PGLn  = X' \dasharrow X$. Choose $x' \in X'$,
  so that $f$ is defined at $x'$ and set $x=f(x')$.  Denote
  by $S$ and $S'$ the stabilizers in $\GLm\times\PGLn$ of $x$ and
  $x'$, respectively.
  Note that $S' \subseteq S \subseteq \GL_m \times \PGLn$.  Since $\GL_m$
  acts transitively on $X$, the projection of $S$ to $\PGLn$ is all of
  $\PGLn$.  On the other hand, we clearly have $S' \subseteq G \times \{ 1 \}$. 
  Consequently, $\dim(S) - \dim(S') \ge \dim(\PGLn) = n^2 -1$, 
  and if $O'$ is a $\GL_m \times \PGLn$-orbit in general 
  position in $X'$, then $\dim(O') - \dim(X) \ge n^2 -1$. 
  We thus conclude that
  \begin{align*}
    \trdeg_{\Z(A)} L &= \trdeg_k \, L - \trdeg_k \Z(A)\\
       &= \dim(X'/\PGLn) - \dim(X/\PGLn) \\
       &=\dim(X') - \dim(X) \ge \dim(O') - \dim(X) \ge n^2 -1 \, ,
  \end{align*}
  as claimed.
\end{proof}

\begin{remark} \label{rem.split-b}
  One can show directly that the $\GL_m$-splitting field $L$ for $A =
  \UDmn$ given by Proposition~\ref{prop3} satisfies the inequality
  of Proposition~\ref{prop.split}(b) (assuming, of course, that $m \ge
  n^2$). Indeed, since $G$ has a dense orbit in $W = X/\PGLn$, for $w
  \in W$ in general position,
  \[ \dim \, \Stab_{\GL_m}(w) = \dim(\GL_m) - \dim(W) 
                              = m^2 - \dim(X/\PGLn)\,. \]
  Since the associated variety $X = \Mnm$ has dimension $mn^2$, this
  yields
  \[ \trdeg_{\Z(A)} \, L = \dim \, \Stab_{\GL_m}(w) 
                         = m(m-n^2) + (n^2 - 1) \ge n^2 - 1 \, , \]
  as claimed.
\end{remark}

\section{An example: the $\GL_2$-action on $\UD(2,2)$}
\label{sect.ex2x2}

In this section we will use Theorem~\ref{thm2} to study the natural
$\GL_m$-action on the universal division algebra $\UDmn$, described in 
Example~\ref{ex.universal}, for $m = n = 2$.
Note that this case exhibits some special features that do not
recur for other values of $m$ and $n\geq 2$; see 
Proposition~\ref{prop.split}(a) (for $m \ge n^2$) and
\cite{rv5} (for $m \le n^2 -1$). 

\begin{prop} \label{prop.UD(2,2)}
  The fixed algebra $\UD(2, 2)^{GL_2}$ is a non-central subfield of
  $\UD(2, 2)$ of transcendence degree~$1$ over $k$.
\end{prop}

Recall from Example~\ref{ex.universal} that the $\GL_2$-action on
$\UD(2, 2)$ is defined as follows.  Denote by $X$ and $Y$ the two
generic $2\times 2$ matrices generating $\UD(2,2)$.  Then for $g \in
\GL_2$, we have $g(X)=\alpha X + \beta Y$, and $g(Y)= \gamma X +
\delta Y$, where
$g\inv=\bigl(\begin{smallmatrix}\alpha&\beta\\
  \gamma&\delta\end{smallmatrix}\bigr)$.
Recall also that the associated variety for the $\GL_2$-action on
$\UD(2,2)$ is $X=(\M_2)^2$.  In order to use Theorem~\ref{thm2}
to prove Proposition~\ref{prop.UD(2,2)}, we first need to determine
the stabilizer in general position for the $\GL_2\times\PGL_2$-action
on $(\M_2)^2$.

\begin{lem} \label{lem.UD(2,2)}
  For $x \in (\M_2)^2$ in general position, $\Stab_{\GL_2 \times
    \PGL_2}(x)$ is isomorphic to $\bbZ/ 2 \bbZ$.
\end{lem}

\begin{proof}
  By Lemma~\ref{lem7.1}, $\Stab_{\GL_2 \times \PGL_2}(x)$ is
  isomorphic to $\Stab_{\GL_2}(y)$ for the $\GL_2$-action on $W =
  X/\PGL_2$, which is a birational model for the $\GL_2$-action on the
  center $Z$ of $\UD(2, 2)$. In this case there is a
  particularly simple birational model, which we now describe.
  
  It is well known that $Z$ is freely generated (as a field
  extension of~$k$) by the five elements $\tr(X)$, $\tr(Y)$,
  $\tr(X^2)$, $\tr(Y^2)$ and $\tr(XY)$; see~\cite[Theorem
  2.2]{procesi1}.  In other words, the categorical (and, hence, the
  rational) quotient for the $\PGL_2$-action is $\bbA^5$. The group
  $\GL_2$ acts on $\bbA^5$ linearly. In fact, the representation of
  $\GL_2$ on $\bbA^5 = X/\PGL_2$ can be decomposed as $V_2 \oplus
  V_3$, where $V_2$ is the natural $2$-dimensional representation (we
  can think of it as $\Span_k(\tr(X), \tr(Y))$) and $V_3$ is its
  symmetric square. (We can think of $V_3$ as $\Span_k(\tr(X^2),
  \tr(Y^2), \tr(XY))$.)
  
  The question we are asking now reduces to the following: What is the
  stabilizer, in $\GL_2$, of a pair $(v, q)$, in general position,
  where $v$ is a vector in $k^2$ and $q$ is a quadratic form in $2$
  variables?  Indeed, since $\GL_2$ acts transitively on
  non-degenerate quadratic forms in two variables, we may assume that
  $q$ is a fixed form of rank $2$, e.g., $q = x^2 + y^2$.  The
  stabilizer of $q$ is thus the orthogonal group $O_2$, and our
  question further reduces to the following: what is the stabilizer in
  general position for the natural linear action of $O_2$ on $k^2$?
  The answer is easily seen to be $\bbZ/2 \bbZ$, where the non-trivial
  element of $\Stab_{O_2}(v)$ is the orthogonal reflection in $v$; see
  Example~\ref{ex.O_2}.
\end{proof}

\begin{proof}[Proof of Proposition~\ref{prop.UD(2,2)}]
  Note that the $\GL_2$-action on $X=(\M_2)^2$ is generically free (it
  is isomorphic to the direct sum of $4$ copies of the natural
  $2$-dimensional representation of $\GL_2$). Thus the image of the
  stabilizer $\Stab_{\GL_2 \times \PGL_2}(x)$ under the natural
  projection to the second factor is $\bbZ/ 2 \bbZ$.  Since this image
  is non-trivial, Theorem~\ref{thm2}(a) tells us that $\UD(2,
  2)^{\GL_2}$ is not a division subalgebra of $\UD(2,2)$ of degree
  $2$. In other words, it is a subfield of $\UD(2, 2)$. On the other
  hand, Theorem~\ref{thm2}(b) tells us that $\UD(2, 2)^{\GL_2}$ is not
  contained in the center~$Z$ of $\UD(2, 2)$.  Indeed, every subgroup
  of $\PGL_2$ of order $2$ is contained in a torus.  Hence,
  $\Stab_{\GL_2 \times \PGL_2}(x)$ is contained in $\GL_2 \times T_x$,
  where $T_x$ is a maximal torus of $\PGL_2$.  It follows from
  Theorem~\ref{thm2}(b) that the subfield $\UD(2, 2)^{\GL_2}$ is not
  central in $\UD(2,2)$.
  
  Finally, note that $\UD(2,2)^{\GL_2}$ is algebraic over $Z^{\GL_2}$,
  since the minimal polynomial of any element of $\UD(2,2)^{\GL_2}$
  over $Z$ is unique, so must have coefficients in $Z^{\GL_2}$.  It
  follows from Lemmas~\ref{lem7.1} and~\ref{lem.UD(2,2)} that the
  $\GL_2$-action on $W=X/\PGL_2$ has a finite stabilizer in general
  position.  Hence the transcendence degree of $Z^{\GL_2}=
  k(X/\PGL_2)^{\GL_2}$ (over $k$) is $\dim(X/\PGL_2)-\dim(\GL_2)=1$.
\end{proof}

\begin{remark} \label{rem-ex1a}
  This argument also shows that $\UD(2, 2)^{\SL_2}$ is a division
  algebra of degree $2$.
\end{remark}

\begin{remark} \label{rem-ex1b}
  One can exhibit an explicit non-central $\GL_2$-fixed element of
  $\UD(2, 2)$ as follows.  Let
  \[ S_3(A_1, A_2, A_3) = \sum_{\sigma \in \Sym_3} (-1)^\sigma
   A_{\sigma(1)} A_{\sigma(2)} A_{\sigma(3)}
  \]
  be the standard polynomial in three variables; cf.,~\cite[p.
  8]{rowen:PI}.  Set $a = [X, Y]=XY-YX$ and $b = S_3(X, Y, a)$. Using
  the fact that $[A_1,A_2]$ and $S_3(A_1,A_2,A_3)$ are multilinear and
  alternating, it is easy to see that for $g\in\GL_2$, $g(a) =
  a/\!\det(g)$ and $g(b) = b/\!\det^2(g)$.  Specializing $X$ to
  $\bigl(\begin{smallmatrix}1&0\\0&0\end{smallmatrix}\bigr)$ and $Y$
  to $\bigl(\begin{smallmatrix}1&1\\1&0\end{smallmatrix}\bigr)$, an
  elementary computation shows that $a$ and $b$ specialize to
  $\bigl(\begin{smallmatrix}\hphantom{-}
    0&1\\-1&0\end{smallmatrix}\bigr)$ and
  $\bigl(\begin{smallmatrix}-4&\hphantom{-}0\\
    \hphantom{-}0&-2\end{smallmatrix}\bigr)$, respectively. This shows
  that $\det(a)\neq0$ and that $b$ is non-central. Now, $b/\!\det(a)$
  is a non-central $\GL_2$-fixed element of $\UD(2,2)$.
  
  Note also that $a$ and $b$ are non-commuting $\SL_2$-invariant
  elements of $\UD(2, 2)$. This gives an explicit proof of
  Remark~\ref{rem-ex1a}.
\end{remark}

\section{An example: a finite group action on a cyclic algebra}
\label{sect.ex-mult}

In this section we present an example of a finite group action 
on a cyclic algebra. This example illustrates Lemma~\ref{lem7.1} 
and Theorem~\ref{thm4} and, in particular, shows that 
the converse to Theorem~\ref{thm4}(b) is false.

Let $p$ be a prime integer, and $\zeta$ a primitive $p$-th root of
unity in $k$.  Let $P=k\{x,y\}$ be the skew-polynomial ring with
generators $x$ and $y$, subject to the relation 
\[xy=\zeta yx\,.\]
Let $A$ be the division algebra of fractions of $P$; it is a central
simple algebra of degree~$n=p$.  Note that $A$ is the symbol algebra
$(u,v)_p$ whose center is $\Z(A) = k(u,v)$, where $u = x^p$ and $v =
y^p$ are algebraically independent over $k$.

For $(a,b)\in(\Zp)^2$, define an automorphism $\sigma_{(a,b)}$ of $A$
by
\begin{equation}\label{e.action-of-K}
  \sigma_{(a,b)}(x)=\zeta^a x \quad\text{and}\quad
  \sigma_{(a,b)}(y)=\zeta^b y\,.
\end{equation}
These automorphisms of $A$ form a group~$K$ which is isomorphic to
$(\Zp)^2$.  Next, we define an automorphism $\tau$ of $A$ by
\begin{equation}\label{e.action-of-tau}
  \tau(x)=y \quad\text{and}\quad \tau(y)=x\inv y\inv\,.
\end{equation}
Note that $\tau$ is well-defined since
\[\tau(x)\tau(y)-\zeta\tau(y)\tau(x)
=y(yx)\inv-\zeta x\inv=y(\zeta\inv xy)\inv-\zeta x\inv=0\,.\]
Elementary calculations show that $\tau$ has order three, and that
$\tau\inv\sigma_{(a,b)}\tau=\sigma_{(b,-a-b)}$.  Consequently, the
subgroup $G$ of automorphisms of $A$ generated by $K$ and $\tau$ is a
semidirect product $G=K\sdp H$, where $K\simeq(\Zp)^2$ and
$H=\langle\tau\rangle\simeq\bbZ/3\bbZ$.

One easily checks that sending $\tau$ to the matrix
$\bigl(\begin{smallmatrix}0&-1\\1&-1\end{smallmatrix}\bigr)$ defines a
representation
\[\phi_p\colon H\to\SL_2(\Zp)\,,
\]
and thus an action of $H$ on $(\Zp)^2$.

Let $X$ be the $G \times \PGLn$-variety associated to
the action of $G$ on the central simple algebra $A$ of degree $n=p$.
That is, $X$ is an irreducible $G\times\PGLn$-variety which is
$\PGLn$-generically free, and $A$ is $G$-equivariantly isomorphic to
$k_n(X)$.

\begin{prop} \label{prop.example1}
  \begin{itemize}
  \item[\textup{(a)}]For $x \in X$ in general position, there exists a
    maximal torus $T_x$ of $\PGLn$ such that $\Stab_{G \times \PGLn}
    (x) \subseteq G \times \N(T_x)$.
  \item[\textup{(b)}]$A$ has a $G$-invariant maximal subfield if and
    only if the 2-dimensional representation $\phi_p\colon
    H\to\SL_2(\Zp)$ is reducible over $\Zp$.
  \item[\textup{(c)}]The converse to Theorem~\ref{thm4}(b) is false.
  \end{itemize}
\end{prop}

Before we proceed with the proof, two remarks are in order.  First of
all, every finite group action on a central simple algebra is automatically
geometric (and algebraic).

Secondly, an explicit model for $X$ is not immediately transparent (a
description of $X$ as a $\PGLn$-variety can be found in~\cite[Lemma
5.2]{ry}).  On the other hand, the $G$-variety $W$ associated to the
$G$-action on the center of $A$ (see the beginning of
Section~\ref{sect.center}) is easy to describe: We can take $W$ to be
the two-dimensional torus $W=(k^*)^2=\Spec(k[u,v,u\inv,v\inv])$, where
as before, $u=x^p$ and $v=y^p$. It follows from~\eqref{e.action-of-K}
and~\eqref{e.action-of-tau} that the $K$-action on $W$ is trivial, and
that the action of $\tau$ is induced from $\tau(u)=v$, $\tau(v)=(x\inv
y\inv)^p=\epsilon\cdot u\inv v\inv$, where $\epsilon=1$ if $p>2$ and
$\epsilon=-1$ if $p=2$.

\smallskip
We now proceed with the proof of Proposition~\ref{prop.example1}.

\begin{proof} (a) Since $G$ is a finite group,
  $\Stab_G(w)$, for $w \in W$ in general position, is precisely the
  kernel of the $G$-action on $W$.  We claim that the kernel is equal
  to $K$.  That it contains $K$ is immediate
  from~\eqref{e.action-of-K}, since every element of $K$ preserves
  both $u = x^p$ and $v = y^p$.  On the other hand, the
  $H$-action on $W$ is faithful, because $H$ is a simple
  group acting nontrivially on $\Z(A)=k(W)$.  We have thus shown that
  $\Stab_G(w) = K$ for $w \in W$ in general position.

  By Lemma~\ref{lem7.1}, $\Stab_{G \times \PGLn}(x) \simeq K \simeq
  (\bbZ/ p \bbZ)^2$ for $x$ in general position in $X$. In particular,
  the projection of this group to $\PGLn$ is a finite abelian
  subgroup of $\PGLn$.  By~\cite[II.5.17]{ss}, every finite abelian
  subgroup of $\PGLn$ lies in the normalizer of a maximal torus $T_x$.
  Thus
  \[ \Stab_{G \times \PGLn}(x) \subseteq G \times \N(T_x) \, ,
  \]
  as claimed.

  \smallskip (b) First we will describe the $K$-invariant maximal
  subfields of $A$, then determine which ones of them are also
  invariant under $H$.  Note that since $A$ is a division algebra of
  prime degree $p$, every nontrivial field extension $L$ of the center
  $\Z(A)$ is a maximal subfield of $A$.

  The group $K\simeq(\bbZ/p \bbZ)^2$ acts trivially on $\Z(A)$; its
  action on $A$ decomposes as a direct sum of $p^2$ one-dimensional
  character spaces $\Span_{\Z(A)}(x^i y^j)$, where $0 \le i, j \le
  p-1$.  These spaces are associated to the $p^2$ distinct characters
  of $(\bbZ/p \bbZ)^2$; hence, every $K$-invariant $\Z(A)$-vector
  subspace $L$ contains $x^i y^j$ for some $0 \le i, j \le p-1$.
  Moreover, if $L$ is a $K$-invariant maximal subfield of $A$ then
  $\Z(A)(x^i y^j) \subseteq L$, where $0 \le i, j \le p-1$ and $(i, j)
  \ne (0, 0)$.  Since $[L:\Z(A)] = p$ and $x^i y^j \not \in \Z(A)$, we
  conclude that $L = \Z(A)(x^i y^j)$. We will denote 
  $\Z(A)(x^i y^j)$ by $L_{(i,j)}$.

  Now suppose $(i, j)$ and $(r, s)$ are non-zero elements of $(\bbZ/p
  \bbZ)^2$.  We claim that $L_{(i,j)} = L_{(r,s)}$ if and only if
  $(i, j)$ and $(r, s)$ are proportional, i.e., if and only if they
  lie in the same 1-dimensional $\bbZ/ p \bbZ$-subspace of $(\bbZ/ p
  \bbZ)^2$. Indeed, if $(i, j)$ and $(r, s)$ are proportional then up
  to a multiple from $\Z(A)$, $x^i y^j$ and $x^r y^s$ are
  powers of one another. Since neither one is central, they generate
  the same maximal subfield.
  Conversely, since a maximal subfield has dimension $p$ over $\Z(A)$, 
  it can contain only $p-1$ distinct $x^iy^j$ with
  $(0,0)\neq(i,j)\in(\Zp)^2$.  Since there are $p-1$ nonzero
  $\Zp$-multiples of $(i,j)$, this proves the claim.

  We have thus shown that the $K$-invariant maximal subfields of $A$
  are in bijective correspondence with 1-dimensional
  $\bbZ/p\bbZ$-subspaces of $(\bbZ/ p \bbZ)^2$: a 1-dimensional subspace
  $V$ corresponds to the maximal subfield $L_V = \Z(A)(x^i y^j)$,
  where $(i, j)$ is a non-zero element of $V$.

  It is clear from~\eqref{e.action-of-tau} that $\tau(L_{V}) =
  L_{\tau(V)}$, where $\tau$ acts on $(\Zp)^2$ via the
  representation~$\phi_p$.  To sum up: $A$ has a maximal $G$-invariant
  subfield $\Longleftrightarrow$ $\tau$ preserves one of the
  $K$-invariant maximal subfields $L_{V}$ $\Longleftrightarrow$
  $(\bbZ/p \bbZ)^2$ has a $\tau$-invariant 1-dimensional $\bbZ/p
  \bbZ$-subspace $V \subset (\bbZ/p \bbZ)^2$ $\Longleftrightarrow$ the
  representation $\phi_p$ of $H$ is reducible.

  \smallskip (c) In view of parts (a) and (b) it suffices to show
  that the representation $\phi_p$ of $H$ is irreducible if and only
  if $p \equiv 2 \pmod{3}$.  If $p=3$, $\phi_p$ is reducible, since
  in this case 
  $\bigl(\begin{smallmatrix}\hphantom{-}1\\-1\end{smallmatrix}\bigr)$ is an
  eigenvector for the matrix
  $\bigl(\begin{smallmatrix}0&-1\\1&-1\end{smallmatrix}\bigr)$.
  Now assume that $p\neq3$.  Then Maschke's theorem implies 
  that $\phi_p(\tau)$ is
  diagonalizable over the algebraic closure of $\Zp$.  The eigenvalues
  of $\phi_p(\tau)$ are then necessarily third roots of unity,
  including at least one primitive third root of unity.

  Thus the action of $H$ on $(\Zp)^2$ is reducible
  $\Longleftrightarrow$ $\phi_p(\tau)$ is diagonalizable over $\Zp$
  $\Longleftrightarrow$ the eigenvalues of $\phi_p(\tau)$ belong to $\Zp$
  $\Longleftrightarrow$ $\Zp$ contains a primitive third root of unity
  $\Longleftrightarrow$ $3\mid p-1$.

  Consequently, the representation $\phi_p$ irreducible if and only if
  $p\equiv 2 \pmod{3}$.
\end{proof}

\appendix
\section{Inner actions on division algebras}
\label{sect.inner}

In this appendix we continue to assume that $k$ is an algebraically
closed base field of characteristic zero, and that every division algebra
is finite-dimensional over its center, which in turn is a finitely
generated field extension of~$k$.  (Some of the lemmas below hold in
greater generality; see Remark~\ref{rem.inner}.)  Our main
result is the following theorem.

\begin{thm} \label{thm.inner.2}
  Let $G$ be an algebraic group acting on a division algebra $D$ of
  degree~$n$ by inner automorphisms.  Then the kernel $N$ of this
  action contains the connected component $G^0$ of $G$, and $G/N$ is a
  finite abelian $n$-torsion group.
\end{thm}

Here the algebraic group $G$ is treated as an abstract group;
in particular, the (inner) action of $G$ on $D$ is not
assumed to be algebraic or geometric. Consequently, 
our proof has a rather different flavor from the other arguments
in this paper. Instead of using algebraic geometry, we exploit,
in the spirit of \cite{ta}, the fact that connected 
algebraic groups are generated, as abstract groups, by their
divisible subgroups. Note that the special case 
of Theorem~\ref{thm.inner.2}, where $G$ is a torus 
is proved in~\cite[Corollary 5.6]{ta}.  

Before we prove Theorem~\ref{thm.inner.2}, we deduce an easy
consequence.

\begin{cor} \label{cor.inner.3}
 Let $G$ be an algebraic group acting faithfully and geometrically on
  a division algebra $D$ of degree~$n$.  Then the normal subgroup of~$G$
  acting by inner automorphisms is a finite abelian $n$-torsion
  group.
 \end{cor}

\begin{proof}
  Since $G$ acts geometrically, the normal subgroup~$H$ of $G$
  consisting of the elements acting by inner automorphisms (i.e.,
  acting trivially on~$\Z(D)$) is closed, so itself an algebraic
  group.  Now apply Theorem~\ref{thm.inner.2} to the
  faithful action of $H$ on $D$.
\end{proof}

We now turn to the proof of Theorem~\ref{thm.inner.2}, beginning with
two lemmas.

\begin{lem} \label{lem.divisible}
  The group of inner automorphisms of a division algebra contains no
  divisible subgroups.
\end{lem}

\begin{proof} 
  Assume to the contrary that there is a nontrivial divisible group
  $H$ acting faithfully on a division algebra $D$ by inner
  automorphisms.  By \cite[Corollary~3.2]{ta}, the torsion subgroup of
  $H$ acts trivially on $D$, so it must be trivial.  Hence $H$ is a
  torsion-free divisible group, i.e., a direct sum of copies of
  $(\bbQ,+)$; cf.~\cite[5.2.7]{scott}.  By \cite[Lemma 3.3(a)]{ta},
  there is a subfield $L$ of $D$ containing the center $K$ of $D$ such
  that $H$ embeds into $L^*/K^*$.  Thus $(\bbQ,+)$ embeds into
  $L^*/K^*$.  By~\cite[Lemma 5.5]{ta}
  \footnote{%
    We take the opportunity to correct an error in the proof of
    \cite[Lemma 5.5]{ta}.  The third paragraph of that proof should
    read: ``If $\pi_i\circ\phi$ is not injective, its image is a
    torsion group.  Since $\pi\circ\phi$ is injective,
    $\pi_i(\phi(\bbQ))$ is not torsion for some $i$.  Hence, for this
    $i$, $\psi=\pi_i\circ\phi$ is injective, so that $\psi(\bbQ)$ is
    nontrivial.  Thus by the argument in the previous paragraph,
    $\psi(\bbQ)$ is not contained in $K^*$."}, this implies that $K$
  is not finitely generated over the algebraically closed field~$k$,
  a contradiction.
\end{proof}

\begin{lem} \label{lem.inner}
  Let $D$ be a division algebra of degree $n$ whose center $K$
  contains all roots of unity.
  \begin{itemize}
  \item[\textup(a)]Suppose $x \in D$ has the following properties:
    $\det(x) = 1$, and $x^m \in K$ for some integer $m \ge 1$. Then $x
    \in K$.
  \item[\textup(b)]If $G$ is a finite group acting faithfully on $D$
    by conjugation, then $G$ is an abelian $n$-torsion group.
  \end{itemize}
\end{lem}

As the statement of the lemma implies, here $K$ is not assumed to
contain an algebraically closed base field.

\begin{proof}
  (a) Suppose $x^m = a$ for some integer $a \in K$. Taking the
  determinant (i.e., reduced norm) on both sides, we obtain $a^n = 1$.
  Thus, after replacing $m$ by $mn$, we may assume $x^m = 1$. Since
  the polynomial $f(t) = t^m - 1$ splits over $K$, we conclude that $x
  \in K$.

  \smallskip
  
  (b) Suppose $g \in G$ acts by conjugation by $d_g$. Then for every
  $g, h \in G$, the commutator $x = d_g d_h d_g^{-1} d_h^{-1}$
  satisfies the conditions of part (a), where $m$ can be taken to be
  the order of $ghg^{-1}h^{-1}$ in $G$. Thus $x \in K$ and
  consequently, $g$ and $h$ commute in $G$. This shows that $G$ is
  abelian.
  
  To prove that $G$ is $n$-torsion, choose $g \in G$ and consider the
  element $x = (d_g)^n/\det(d_g)$.  Once again, $x$ satisfies the
  conditions of part (a), with $m$ the order of $g^n$ in $G$. Thus $x
  \in K$, and consequently, $g^n = 1$ in $G$, as claimed.
\end{proof}

\begin{proof}[Proof of Theorem~\ref{thm.inner.2}]
  Let $S$ be a torus of $G$, or a closed subgroup which is isomorphic
  to $(k,+)$.  We claim that $S \subseteq N$.  Since $S$ is a divisible
  group, so is $S/N \cap S$; cf.~\cite[5.2.19]{scott}.  Since $S/N
  \cap S$ acts faithfully on $D$, Lemma~\ref{lem.divisible} tells us
  that $S/N \cap S = \{ 1 \}$, i.e., $S \subseteq N$, as claimed.
  
  Now recall that every element $g \in G^0$ has a Jordan decomposition
  product $g = g_s g_u$, where $g_s$ is semisimple and $g_u$ is
  unipotent; cf., e.g.,~\cite[Theorem 15.3]{humphreys}.  Since $g_s$
  lies in a torus of $G$, $g_s \in N$. Similarly, $g_u \in N$; cf.,
  e.g.,~\cite[Lemma~15.1C]{humphreys}.  Thus $G^0 \subseteq N$, as
  claimed.  The desired conclusion now follows from
  Lemma~\ref{lem.inner}.
\end{proof}

\begin{remark}\label{rem.inner}
  Lemmas \ref{lem.divisible} and \ref{lem.inner} also hold in prime
  characteristic, and so does Theorem~\ref{thm.inner.2}, provided
  $G$ is reductive (since then $G^0$ is generated as abstract group by
  the tori it contains).
\end{remark}
 
\section{Regular actions on prime PI-algebras}
\label{sect.prime.PI}

It is a consequence of Posner's theorem that every prime
PI-algebra~$R$ of PI-degree~$n$ can be realized as a subalgebra of
$n\times n$-matrices over some commutative domain~$C$.  Given an
action of a group~$G$ on~$R$, it is natural to ask whether
one can always find such an embedding $R\hookrightarrow\Mn(C)$ which
is $G$-equivariant for some action of $G$ on $\Mn(C)$.  We now deduce
from Theorem~\ref{thm3} a rather strong affirmative answer in the case
of regular actions of algebraic groups (see Definition~\ref{def.regular})
on prime PI-algebras.  Such actions were extensively 
studied in \cite{vonessen:memoir} and \cite{vonessen:trans}.  

\begin{prop}\label{prop:regular}
  Let $R$ be a prime PI-algebra of PI-degree~$n$, which is finitely
  generated as $k$-algebra.  Let $G$ be an algebraic group
  acting regularly on~$R$.  Then there is a finitely generated
  commutative $k$-algebra $C$ which is a domain, and a regular action
  of $G$ on $C$ such that $R$ embeds $G$-equivariantly into
  $\Mn\otimes_kC$.  Here $G$ acts trivially on $\Mn$.
\end{prop}

In the case where $G$ is a torus, this assertion was proved in
\cite[Corollary~9]{vonessen:solv}. 

\begin{proof}
  Let $A$ be the total ring of fractions of $R$; it is a central
  simple algebra of degree~$n$, and $G$ acts algebraically on $A$.
  Note that since $R$ is finitely generated as $k$-algebra, the center
  of $A$ is a finitely generated field extension of~$k$. By
  Theorem~\ref{thm3}, there is a $G$-splitting field $L = k(X_0)$ 
  for $A$, where $X_0$ is an affine $G$-variety, i.e., the
  $G$-action on $L$ is algebraic; cf.\ Definition~\ref{defn:algebraic-action}.
  This gives rise to a $G$-equivariant embedding $\varphi\colon
  R\to\Mn\otimes_k L=A'$.  Hence $G$ also acts algebraically on $A'$,
  so that $A'$ contains a unique largest subalgebra $S_{A'}$ on which
  $G$ acts regularly, and which contains every subalgebra of $A'$ on
  which $G$ acts regularly.  Denote by $S_L$ the corresponding
  subalgebra of~$L$.  Since $S_{A'}$ contains $\Mn\otimes_kk$, it
  follows that $S_{A'}=\Mn\otimes_kS_L$.  Since $G$ acts regularly on
  $\varphi(R)$, $\varphi(R)\subseteq\Mn\otimes_kS_L$.  Since $R$ is
  finitely generated, and since $G$ acts regularly on $S_L$, there is
  a finitely generated $G$-invariant subalgebra $C$ of $S_L$ such that
  $\varphi(R)\subseteq\Mn\otimes_kC$.
\end{proof}                 
             

\end{document}